\documentclass[a4paper,fleqn]{cas-sc}
\usepackage[utf8]{inputenc} 
\usepackage[authoryear]{natbib}

\usepackage{textcmds}

\usepackage{graphicx} 
\usepackage{tikz}
\usetikzlibrary{positioning}
\usepackage{booktabs}
\usepackage{supertabular}
\usepackage{longtable}
\usepackage{caption}
\usepackage{subcaption}
\usepackage[toc,page]{appendix}
\usepackage{threeparttable}
\usepackage{multirow}
\usepackage{tabularx}
\usepackage{chngcntr}
\usepackage{appendix}
\newcolumntype{Y}{>{\centering\arraybackslash}X}

\usepackage{booktabs} 
\usepackage{array} 
\usepackage{paralist} 
\usepackage{verbatim} 
\usepackage{amsmath, amssymb}
\usepackage{stix}
\usepackage{pgf}
\usepackage{hyperref}

\definecolor{gray1}{gray}{0.4}
\definecolor{gray2}{gray}{0.7}
\definecolor{gray3}{gray}{0.9}
\definecolor{gray20P}{gray}{0.2}
\definecolor{gray25P}{gray}{0.25}
\definecolor{gray30P}{gray}{0.3}
\definecolor{gray40P}{gray}{0.4}
\definecolor{gray50P}{gray}{0.5}
\definecolor{gray60P}{gray}{0.6}
\definecolor{gray70P}{gray}{0.7}
\definecolor{gray80P}{gray}{0.8}
\definecolor{gray85P}{gray}{0.85}
\definecolor{gray90P}{gray}{0.9}

\newcolumntype{L}{>{$}l<{$}}
\def\variant#1#2#3{\makebox[0.2cm][c]{#1}\,|\,\makebox[0.20cm][c]{#2}\,|\,\makebox[0.2cm][c]{#3}}

\begin{document}
\let\WriteBookmarks\relax
\def\floatpagepagefraction{1}
\def\textpagefraction{.001}

\title[mode=title]{The vehicle routing problem with synchronization constraints and support vehicle-dependent service times}
\shorttitle{VRP with synchronization constraints and support vehicle-dependent service times}    

\shortauthors{Wittwer and Tamke}  
\author[1]{David Wittwer} 
\cormark[1]
\ead{david.wittwer@tu-dresden.de}
\affiliation[1]{organization={TU Dresden, Faculty of Mechanical Science and Engineering},
            city={Dresden},
            postcode={01062}, 
            state={Saxony},
            country={Germany}}

\author[2]{Felix Tamke}\ead{felix.tamke@tu-dresden.de}
\affiliation[1]{organization={TU Dresden, Faculty of Business and Economics},
            city={Dresden},
            postcode={01062}, 
            state={Saxony},
            country={Germany}}
\cortext[cor1]{Corresponding author}
\begin{abstract}
Many production processes require the cooperation of various resources. Especially when using expensive machines, their utilization plays a decisive role in efficient production. In agricultural production or civil construction processes, e.g., harvesting or road building, the machines are typically mobile, and synchronization of different machine types is required to perform operations. In addition, the productivity of one type often depends on the availability of another type.
In this paper, we consider two types of vehicles, called primary and support vehicles. Primary vehicles perform operations and are assisted by at least one support vehicle, with more support vehicles resulting in faster service times for primary vehicles. We call this practical problem the vehicle routing and scheduling problem with support vehicle-dependent service times and introduce two mixed-integer linear programming models. The first represents each support vehicle individually with binary decision variables, while the second considers the cumulative flow of support vehicles with integer decision variables. Furthermore, the models are defined on a graph that allows easy transformation into multiple variants. These variants are based on allowing or prohibiting switching support vehicles between primary vehicles and splitting services among primary vehicles. We show in our extensive computational experiments that: i) the integer representation of support vehicles is superior to the binary representation, ii) the benefit of additional vehicles is subject to saturation effects and depends on the ratio of support and primary vehicles, and iii) switching and splitting lead to problems that are more difficult to solve, but also result in better solutions with higher primary vehicle utilization.
\end{abstract}


\maketitle

\section{Introduction} \label{sec:introduction}

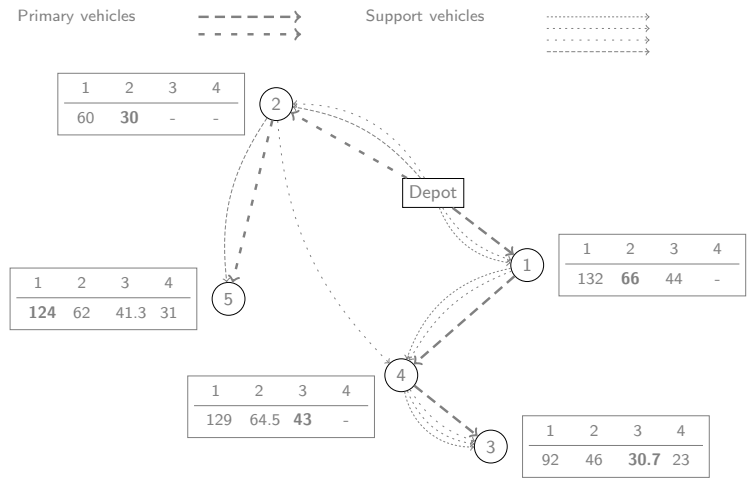
\begin{figure}
    \centering
    \resizebox{0.60\textwidth}{!}
{
\begin{tikzpicture}

    \node[text width=3cm] (primary) at (-2,7.8) {\footnotesize Primary vehicles};
    \node (dummy0_prim_left) at (-0.5, 7.8){};
    \node (dummy1_prim_left) at (-0.5, 7.5){};
    \node (dummy0_prim_right) at (1.5, 7.8){};
    \node (dummy1_prim_right) at (1.5, 7.5){};
    \draw[dash pattern=on 5pt off 3pt, very thick] [->] (dummy0_prim_left) edge node {} (dummy0_prim_right);
    \draw[dash pattern=on 3pt off 5pt, very thick] [->] (dummy1_prim_left) edge node {} (dummy1_prim_right);
    
    \node[text width=3cm] (secondary) at (4,7.8) {\footnotesize Support vehicles};
    \node (dummy0_sec_left) at (5.5, 7.8){};
    \node (dummy1_sec_left) at (5.5, 7.6){};
    \node (dummy2_sec_left) at (5.5, 7.4){};
    \node (dummy3_sec_left) at (5.5, 7.2){};
    
    \node (dummy0_sec_right) at (7.5, 7.8){};
    \node (dummy1_sec_right) at (7.5, 7.6){};
    \node (dummy2_sec_right) at (7.5, 7.4){};
    \node (dummy3_sec_right) at (7.5, 7.2){};
    
    \draw[dash pattern=on 1pt off 1pt, ->] (dummy0_sec_left) edge node {} (dummy0_sec_right);
    \draw [dash pattern=on 1pt off 2pt, ->] (dummy1_sec_left) edge node {} (dummy1_sec_right);
    \draw [dash pattern=on 1pt off 3pt, ->] (dummy2_sec_left) edge node {} (dummy2_sec_right);
    \draw [dash pattern=on 2pt off 1pt, ->] (dummy3_sec_left) edge node {} (dummy3_sec_right);

\node[draw=black] (0) at (3.6587566928488315,4.767929109149463){Depot};
\node[shape=circle,draw=black] (1) at (5.278166920551097,3.52596613168603){1};
\node[shape=circle,draw=black] (2) at (0.955688582726976,6.297792780330559){2};
\node[shape=circle,draw=black] (3) at (4.650874639515099,0.4015031441951322){3};
\node[shape=circle,draw=black] (4) at (3.1087384857087086,1.6295039466773518){4};
\node[shape=circle,draw=black] (5) at (0.13325110272505825,2.9447394749923026){5};

\node[left = 0cm of 5] (times_5) {};
\draw[dash pattern=on 5pt off 3pt, very thick] [->] (0) edge node {}(1);
\draw[dash pattern=on 3pt off 5pt, very thick] [->] (0) edge node {}(2);
\draw[dash pattern=on 5pt off 3pt, very thick] [->] (1) edge node {}(4);
\draw[dash pattern=on 5pt off 3pt, very thick] [->] (4) edge node {}(3);
\draw[dash pattern=on 3pt off 5pt, very thick] [->] (2) edge node {}(5);
\draw [dash pattern=on 1pt off 1pt, thin] [->] (0) edge[bend right=30] node {} (1);
\draw [dash pattern=on 1pt off 1pt, thin] [->] (1) edge[bend right=30] node {} (4);
\draw [dash pattern=on 1pt off 1pt, thin] [->] (4) edge[bend right=40] node {} (3);
\draw [dash pattern=on 1pt off 2pt, thin] [->] (0) edge[bend right=20] node {} (1);
\draw [dash pattern=on 1pt off 2pt, thin] [->] (1) edge[bend right=20] node {} (4);
\draw [dash pattern=on 1pt off 2pt, thin] [->] (4) edge[bend right=30] node {} (3);
\draw [dash pattern=on 1pt off 3pt, thin] [->] (0) edge[bend right=30] node {} (2);
\draw [dash pattern=on 1pt off 3pt, thin] [->] (2) edge[bend right=20] node {} (4);
\draw [dash pattern=on 1pt off 3pt, thin] [->] (4) edge[bend right=20] node {} (3);
\draw [dash pattern=on 2pt off 1pt, thin] [->] (0) edge[bend right=20] node {} (2);
\draw [dash pattern=on 2pt off 1pt, thin] [->] (2) edge[bend right=20] node {} (5);

\node[shape=rectangle,draw, right = 0.25 cm of 1] (times_1){
    \footnotesize
    \begin{tabularx}{3cm}{YYYY}
        1 & 2 & 3 & 4 \\
        \midrule
        132 & \textbf{66} & 44 & -
    \end{tabularx}
};
\node[shape=rectangle,draw, left = 0.25 cm of 2] (times_2){
\footnotesize
    \begin{tabularx}{3cm}{YYYY}
        1 & 2 & 3 & 4 \\
        \midrule
        60 & \textbf{30} & - & -
    \end{tabularx}
};
\node[shape=rectangle,draw, right = 0.25 cm of 3] (times_3){
\footnotesize
    \begin{tabularx}{3cm}{YYYY}
        1 & 2 & 3 & 4 \\
        \midrule
        92 & 46 & \textbf{30.7} & 23
    \end{tabularx}
};
\node[shape=rectangle,draw, below left = -0.2 cm and 0.25 cm of 4] (times_4){
\footnotesize
    \begin{tabularx}{3cm}{YYYY}
        1 & 2 & 3 & 4 \\
        \midrule
        129 & 64.5 & \textbf{43} & -
    \end{tabularx}
};
\node[shape=rectangle,draw, left = 0.25 cm of 5] (times_5){
\footnotesize
    \begin{tabularx}{3cm}{YYYY}
        1 & 2 & 3 & 4 \\
        \midrule
        \textbf{124} & 62 & 41.3 & 31
    \end{tabularx}
};

\end{tikzpicture}
}
    \caption{Solution of the VRP with synchronization constraints and support vehicle-dependent service times for an example with two primary vehicles and four support vehicles. Selected service times at customer locations are printed in bold.}
    \label{fig:Graph_base_model}
\end{figure}

In this paper, we introduce a  vehicle routing problem (VRP) with temporal and spacial synchronized operations of two different types of vehicles. Vehicles of the first type are called primary vehicles. These carry out the actual activities. However, to execute the work, they need the assistance of a second type of vehicle, called support vehicle. Here, the number of support vehicles influences the utilization rate of the primary vehicle and thus the service time. We call this new optimization problem the vehicle routing problem with synchronization constraints and support vehicle-dependent service times. Figure \ref{fig:Graph_base_model} shows an example with five customer requests (circles), two primary vehicles (thick lines), and four support vehicles (thin lines) to illustrate the interaction of primary and support vehicles. The available service times for each customer are given in the tables beside the nodes, and the selected service time is printed in bold. Two support vehicles leave the depot with each primary vehicle to serve customers 1 and 2. After serving customer 2, one support vehicle travels from customer 2 to customer 4 to assist primary vehicle 1 with customers 4 and 3. Primary vehicle 2 travels to customer 5 with the remaining support vehicle. 

\begin{figure}
    \centering
    \includegraphics[width=0.5\textwidth]{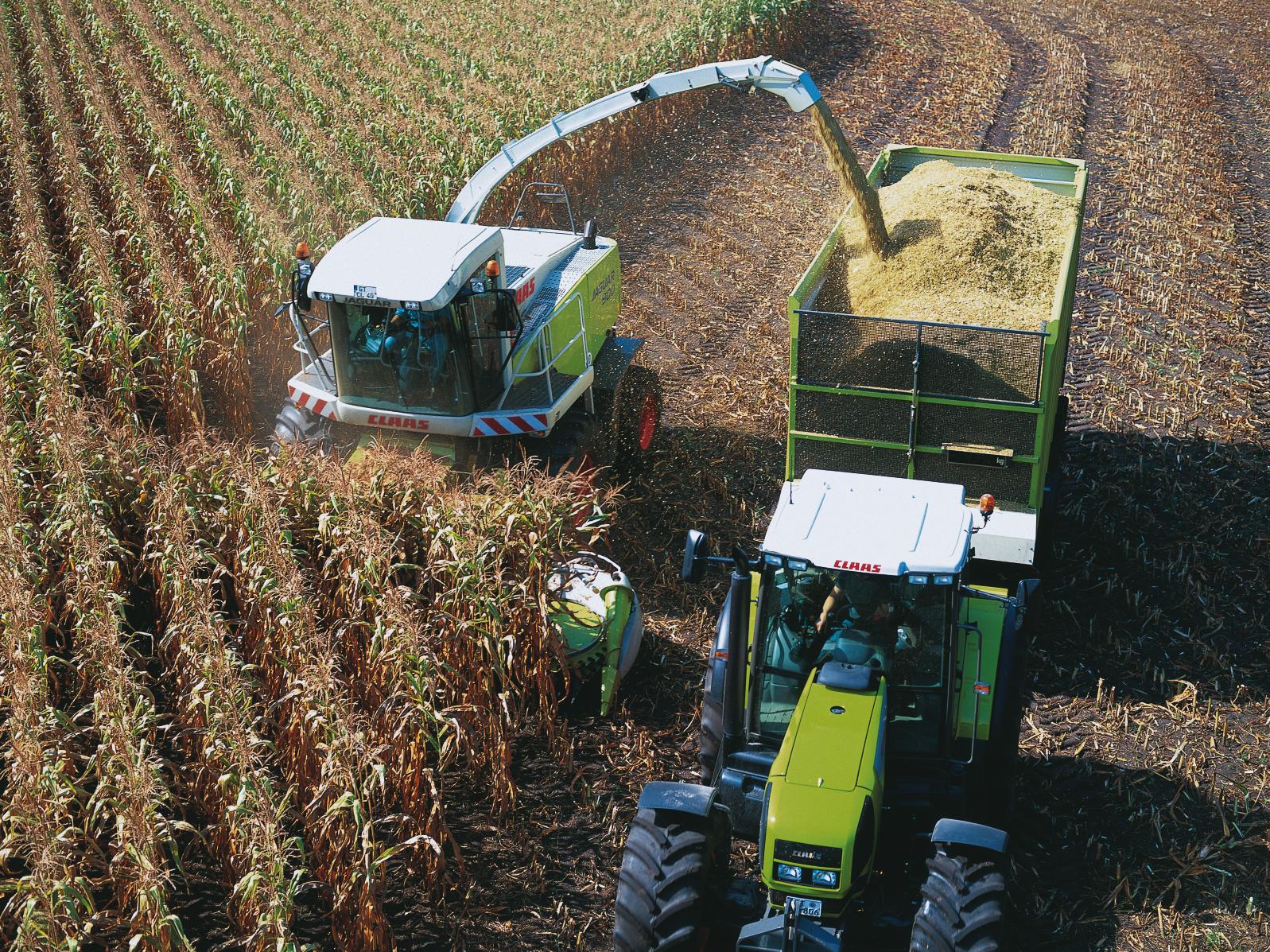}
    \caption{Forage harvester (left) and transporter working together in a cornfield. \citep{gridin_harvester_2005} This file is licensed under the Creative Commons Attribution-Share Alike 3.0 Unported license. To view a copy of this license, visit \url{https://creativecommons.org/licenses/by-sa/3.0/deed.en}.}
    \label{fig:harvester}
\end{figure}

The VRP with synchronization constraints and support vehicle-dependent service times is inspired by the dispatch of forage harvesters and transporters. Figure~\ref{fig:harvester} shows the cooperation between both vehicle types working in a cornfield.
A forage harvester usually cannot store the crops itself. Therefore, transporters as support vehicles are needed to collect the crops during harvest.
They transport the crops to the next storage facility when fully loaded and return to the harvester after unloading the crops. In order to reduce the idle time of the harvester and achieve a high utilization rate, additional transporters are required, since the harvester cannot continue operation without any transporter. Instead of considering the exact movement of transporters during harvest, i.e., to and from the storage facility, we consider that different numbers of transporters result in different service times for the harvester. These service times are customer-specific and are determined in advance, taking into account several factors, such as the distance between the field and the storage facility, the in-field routing of harvesters and transporters, the transporter capacity, and the time for loading and unloading the transporters. Determining the service times is not part of the problem presented here. We refer to studies regarding coverage path planning for harvesters, e.g. \citep{oksanen_coverage_2009, ali_infield_2009}, and path planning for transport vehicles, e.g., \citep{jensen_-field_2012, bochtis_path_2010}.

For the practical forage harvester and transporter dispatching problem two-step heuristics for biomass transfer \citep{amiama_decision_2015, bochtis_vehicle_2010}, heuristics for the calculation of harvester routes \citep{cerdeira-pena_optimised_2017}, and a simulation of the logistical process \citep{amiama_modelling_2015} have been applied. 
The problem presented in this paper is also relevant in other areas, e.g., harvesting with combine harvesters \citep{cordesses_combine_2000}, field fertilizing \citep{hadrich_economic_2010}, concrete casting \citep{yan_effective_2016}, road construction \citep{galic_methodology_2016}, or any other process where a primary resource serves as an enabler (e.g., key, authorization or specially trained staff) and a support resource serves as an accelerator influencing the service time at the customer.
However, to the best of our knowledge, the integrated routing and scheduling of primary vehicles and multiple support vehicles that influence the service time of the primary vehicles have not yet been studied.
We make the following contributions to address this gap:
\begin{enumerate}
    \item We propose two mixed-integer linear programs (MILP) for the synchronized vehicle routing problem with support vehicle-dependent service times. We show that an integer flow representation of support vehicles is superior to a binary flow representation. The underlying graph of our modeling approach allows us to easily transform the model into different variants by simple adjustments.
    \item We provide an analysis of the managerial benefits of primary and support vehicles. Increasing the number of primary and support vehicles leads to reduced completion times. However, the benefit of an additional vehicle decreases due to saturation effects and depends on the ratio of support to primary vehicle count. Thus, if fewer support vehicles are available, an additional primary vehicle is less beneficial than if more support vehicles are available. 
    \item We compare policies that disable/enable switching support vehicles between primary vehicles and splitting service times at customers. We show that switching and splitting lead to problems that are more difficult to solve, but also result in better solutions with higher primary vehicle utilization. In general, enabling the splitting of services is more beneficial than enabling the switching of support vehicles.
\end{enumerate}

The remainder of this paper is structured as follows. We present related literature in Section~\ref{sec:literature}. The synchronized vehicle routing problem with support vehicle-dependent service times and its variants are mathematically described in Section~\ref{sec:problem}. We introduce valid inequalities in Section~\ref{sec:validInequalities} to improve computational performance. Section \ref{sec:ComputationalExperiments} contains computational experiments and discussions. Finally, we summarize the presented work and give an outlook on future work in Section~\ref{sec:Conclusions}.

\section{Related Literature}
\label{sec:literature}
Many relevant optimization problems can be interpreted as either routing or scheduling problems. While routing problems focus on the movement of resources between different locations, scheduling problems emphasize temporal allocation. Given the movement of resources between customer locations, we interpret the novel problem presented here primarily as a variant of the VRP. However, some aspects of the problem, especially different modes to perform activities, are also considered in multi-mode resource-constrained project scheduling problems \citep{coelho_multi-mode_2011, cheng_multi-mode_2015}. For further discussion of the differences and similarities between routing and scheduling problems, we refer the reader to \citep{beck_vehicle_2003}. 

The classical VRP considers a fleet of homogeneous vehicles used to fulfill delivery requests from different customers. Since its introduction several decades ago, it has been expanded in many ways to integrate various practical requirements \citep{toth_vehicle_2014, braekers_vehicle_2016}.
One branch of research includes synchronization requirements between the vehicles concerning different aspects. \citep{drexl_synchronization_2012} identify five different types of synchronization: task, operations, movement, load, and resource synchronization. According to this classification of VRPs with multiple synchronization constraints, the VRP with synchronization constraints and support vehicle-dependent service times contains task synchronization, exact operation synchronization, and resource synchronization. Task synchronization is generally included in any VRP, as delivery tasks or services are distributed among multiple vehicles. Exact operation synchronization occurs, when two or more vehicles are required to coordinate their operations temporally to either simultaneously perform a task or with a certain offset. Here, a primary vehicle and one or more support vehicles begin and end their operations simultaneously and in the same location. Resource synchronization considers vehicles competing for a scarce common resource. In this problem, primary vehicles compete for support vehicles as a mobile resource. Load synchronization is only implicitly taken into account, since the entire load must be transshipped from primary to support vehicles during the operations. There is no movement synchronization because both primary and support vehicles can move independently.

In order to classify the problem presented in this paper accurately, we further differentiate within the group of problems with exact operation synchronization regarding particular characteristics of the resources as well as the customers:
Synchronization may occur between identical resources, between resources of the same type but with another skill set (e.g., worker), or between resources of an entirely different type (e.g., loader and transport vehicle).
If resources can travel between locations by themselves, they are considered active resources. Passive resources, however, require support from other vehicles to move between locations.
A distinctive feature of the customers in this problem is the option to be served in different modes. These modes are determined by the number of resources deployed and influence the customer's service time. This feature is common in scheduling problems but fairly rare in routing problems.
Furthermore, customer demand may be split into sub-jobs, each performed by a different resource or set of resources. If the customer demand exceeds the capacity of the transport vehicles, splitting is mandatory, hence not a deciding factor. If a vehicle can fully satisfy a customer's demand, splitting is optional.

Table \ref{tab:researchGap} provides an overview of the literature reviewed regarding the characteristics above and highlights the research gap. In the following, we describe the problems from the literature in more detail.
\begin{table}
\begin{footnotesize}
\centering
	\begin{threeparttable}
		\begin{tabular}{p{4.3cm}lccccccccl}
	    		\toprule
	    		\multirow{2}{*}[-2em]{Problem classification} & \multirow{2}{*}[-2em]{Publication} & \multicolumn{3}{c}{Operate} & \multicolumn{2}{c}{Move} &
	    		  \multirow{2}{*}[-0.5em]{\rotatebox[origin=c]{90}{Multi-mode}} &
	    		  \multicolumn{2}{c}{Split} & \multirow{2}{*}[-2em]{Method} \\
			    \cmidrule[0.1pt](lr){3-5}
	    		\cmidrule[0.1pt](lr){6-7}
	    		\cmidrule[0.1pt](lr){9-10}
	    		 &  & {\rotatebox[origin=c]{90}{Diff. types}} &
	    		 {\rotatebox[origin=c]{90}{Same type}} &
	    		 {\rotatebox[origin=c]{90}{Hom.}} & {\rotatebox[origin=c]{90}{Active}} & {\rotatebox[origin=c]{90}{Passive}} & & {\rotatebox[origin=c]{90}{Opt.}} & {\rotatebox[origin=c]{90}{Mand.}} &\\
	     		\midrule[0.1pt]
	     		
			
			
			\multirow{2}{3cm}{Manpower allocation} & \citet{li_manpower_2005} && \textbullet && \textbullet &&&&& SA \\ 
			& \citet{dohn_manpower_2009} && \textbullet && \textbullet &&&&& B\&P \\
			\cmidrule[0.1pt](lr){1-2}
			\cmidrule[0.1pt](lr){3-5}
			\cmidrule[0.1pt](lr){6-7}
			\cmidrule[0.1pt](lr){8-8}
			\cmidrule[0.1pt](lr){9-10}
			\cmidrule[0.1pt](lr){11-11}
			
			\multirow{3}{4.3cm}{Home health care routing and scheduling} & \citet{bredstrom_combined_2008} &&  \textbullet  && \textbullet &&&&& MatH\\
			&\citet{mankowska_home_2014} && \textbullet  & & \textbullet &&&&& AVNS \\
			&\citet{lopez-aguilar_linear_2018} && \textbullet &  & \textbullet &&&&& MILP \\
			\cmidrule[0.1pt](lr){1-2}
			\cmidrule[0.1pt](lr){3-5}
			\cmidrule[0.1pt](lr){6-7}
			\cmidrule[0.1pt](lr){8-8}
			\cmidrule[0.1pt](lr){9-10}
			\cmidrule[0.1pt](lr){11-11}
						
			VRP with device synchronization &\citet{cappanera_temporal_2020} & \textbullet & \textbullet && \textbullet & \textbullet &&&& MILP \\
			\cmidrule[0.1pt](lr){1-2}
			\cmidrule[0.1pt](lr){3-5}
			\cmidrule[0.1pt](lr){6-7}
			\cmidrule[0.1pt](lr){8-8}
			\cmidrule[0.1pt](lr){9-10}
			\cmidrule[0.1pt](lr){11-11}
			
			2-Echelon VRP with load transfers & \citet{crainic_models_2009} && \textbullet  & & \textbullet &&&&& Heuristics \\
			\cmidrule[0.1pt](lr){1-2}
			\cmidrule[0.1pt](lr){3-5}
			\cmidrule[0.1pt](lr){6-7}
			\cmidrule[0.1pt](lr){8-8}
			\cmidrule[0.1pt](lr){9-10}
			\cmidrule[0.1pt](lr){11-11}
			
			Concrete delivery & \citet{schmid_hybridization_2010} & \textbullet & \textbullet && \textbullet &&&& \textbullet & VNS, VLNS \\
			\cmidrule[0.1pt](lr){1-2}
			\cmidrule[0.1pt](lr){3-5}
			\cmidrule[0.1pt](lr){6-7}
			\cmidrule[0.1pt](lr){8-8}
			\cmidrule[0.1pt](lr){9-10}
			\cmidrule[0.1pt](lr){11-11}
			
			Military aircraft mission planning & \citet{quttineh_military_2013} && \textbullet && \textbullet &&&&& MILP \\
			\cmidrule[0.1pt](lr){1-2}
			\cmidrule[0.1pt](lr){3-5}
			\cmidrule[0.1pt](lr){6-7}
			\cmidrule[0.1pt](lr){8-8}
			\cmidrule[0.1pt](lr){9-10}
			\cmidrule[0.1pt](lr){11-11}

			\multirow{2}{4cm}{Truck and loader routing} & \citet{rix_column_2015} & \textbullet & \textbullet & & \textbullet & \textbullet &&& \textbullet & CG \\
			&\citet{soares_multiple_2019} & \textbullet & \textbullet &  & \textbullet & \textbullet &&& \textbullet & MatH \\
			\cmidrule[0.1pt](lr){1-2}
			\cmidrule[0.1pt](lr){3-5}
			\cmidrule[0.1pt](lr){6-7}
			\cmidrule[0.1pt](lr){8-8}
			\cmidrule[0.1pt](lr){9-10}
			\cmidrule[0.1pt](lr){11-11}
			
			Gantry and shuttle car scheduling & \citet{fedtke_gantry_2017} &  \textbullet && & \textbullet &&&&&  MatH \\
			\cmidrule[0.1pt](lr){1-2}
			\cmidrule[0.1pt](lr){3-5}
			\cmidrule[0.1pt](lr){6-7}
			\cmidrule[0.1pt](lr){8-8}
			\cmidrule[0.1pt](lr){9-10}
			\cmidrule[0.1pt](lr){11-11}
			
			
			VRP with worker synchronization & \citet{fink_VRPMS_2019} &&& \textbullet && \textbullet & \textbullet &&&  CG \\
			\cmidrule[0.1pt](lr){1-2}
			\cmidrule[0.1pt](lr){3-5}
			\cmidrule[0.1pt](lr){6-7}
			\cmidrule[0.1pt](lr){8-8}
			\cmidrule[0.1pt](lr){9-10}
			\cmidrule[0.1pt](lr){11-11}

			VRP with mobile depots & \citet{hof_intraroute_2021} & \textbullet & & \textbullet & \textbullet &  &&&& ALNS \\
			\cmidrule[0.1pt](lr){1-2}
			\cmidrule[0.1pt](lr){3-5}
			\cmidrule[0.1pt](lr){6-7}
			\cmidrule[0.1pt](lr){8-8}
			\cmidrule[0.1pt](lr){9-10}
			\cmidrule[0.1pt](lr){11-11}

			
			
	     		& This work & \textbullet 	&	& \textbullet			& \textbullet	 & 			& \textbullet & \textbullet &  & MILP \\
	     		\bottomrule
		\end{tabular}
	\end{threeparttable}
 	\caption{Literature review on VRP with exact operation synchronization. Classification according to: operation synchronization between different resource types, resources of the same type but with different abilities, or homogeneous resources, actively or passively traveling resources, multiple serving modes of customers, optional or mandatory splitting of jobs, and solution method applied to the problem.}
	\label{tab:researchGap}
\end{footnotesize}
\end{table}

\citep{li_manpower_2005} present a labor allocation problem where tasks require one or more workers with a different skill set. When more than one worker is required to complete a task, operation synchronization occurs. As teams and vehicle are inseparable, they are considered active resources. The authors apply two construction methods and a simulated annealing (SA) metaheuristic to solve the problem. \citep{dohn_manpower_2009} extend the problem by implementing partial task completion and solve several instances optimally with branch-and-price (B\&P).

A problem frequently discussed in VRP literature is the home health care routing and scheduling problem (HHCRSP). In the variant with operation synchronization, some patients need assistance from two caregivers with certain skills simultaneously. Fixed caregiver vehicle tandems are considered active resources. Publications on this topic consider different properties and different solution methods.
\citep{bredstrom_combined_2008} solve the HHCRSP using a heuristic that iteratively solves restricted MILP problems to improve an initial solution.
\citep{mankowska_home_2014} enable delays at customer nodes but penalize delays and waiting times in the objective function. The authors apply an adaptive variable neighborhood search (AVNS) to solve instances with several hundred customers.
\citep{lopez-aguilar_linear_2018} compare and evaluate different MILP formulations of the problem in terms of their solution efficiency.
\citep{cappanera_temporal_2020} extend the problem by an additional tool that caregivers need for some activities and is passed between them. Using this tool for a task leads to operation synchronization between the tool and the caregivers.  The problem is formulated as MILP with lower bounds and solved with a commercial solver.

\citep{crainic_models_2009} present a 2-echelon VRP in urban logistics. The cargo transfer between support and task vehicles at transfer stations is synchronized. The authors apply a heuristic hierarchical decomposition approach, containing one model each for routing the support and the task vehicles.

\citep{schmid_hybridization_2010} study concrete deliveries from different factories to construction sites. Since demand exceeds vehicle capacity, splitting of the demand is mandatory. At a construction site, only one vehicle may unload at a time. When the vehicle leaves, the next vehicle should already be on site so that filling is continuous across all vehicles. In addition, at some construction sites, support vehicles are required to arrive before the first transporter and leave after the last transporter. The problem is solved with a variable neighborhood search (VNS) and a variable large neighborhood search (VLNS) that improves the initial solution calculated by the VNS. 

\citep{quttineh_military_2013} study a military aircraft mission planning problem. The authors model the problem as a VRP with synchronization constraints, where the target illumination by one aircraft and another aircraft's attack require exact operation synchronization. Small instances of the problem are solved optimally with a commercial solver.

\citep{rix_column_2015} present a log-truck scheduling problem. Transporters and loaders are synchronized at loading points so that one of each (active) vehicle type is present throughout the loading process. Due to high transport demands, splitting is mandatory. The authors solve the problem with a column generation (CG) approach.
\citep{soares_multiple_2019} also study a log-truck scheduling problem. However, loaders cannot move autonomously and are considered passive vehicles. A ``fix and optimize'' matheuristic (MatH) is used to solve the problem. 

\citep{fedtke_gantry_2017} describe a gantry crane scheduling problem that occurs at transshipment points of freight trains. Here, two gantry cranes operating in parallel are deployed. Transports from crane to crane require a  shuttle vehicle that moves between the cranes. This makes synchronization between shuttle vehicles and cranes necessary. The authors solve the problem using a dynamic programming-based matheuristic.

\citep{fink_VRPMS_2019} define a vehicle routing problem with worker and vehicle synchronization to complete various activities at an airport. Workers are passive resources and therefore require vehicles to move between sites, while vehicles cannot travel empty. Workers can perform tasks in different modes: A higher number of workers leads to shorter completion times. By enabling multiple workers to start a task simultaneously, operation synchronization is modeled implicitly. The authors solve the problem by using a column generation-based heuristic approach.

\citep{hof_intraroute_2021} introduce a VRP with mobile depots. Task vehicles can refill their loads and fuel at support vehicles that serve as mobile depots. Therefore, the operations between both vehicles are synchronized. The problem is solved by an Adaptive Large Neighborhood Search (ALNS).

We can state that some aspects of the VRP with synchronization constraints and support vehicle-dependent service times have been studied individually. However, the synchronized routing of two active resource types, primary and support vehicles, where both vehicle types are needed for operations while the availability of one resource affects the service time of the other, has not yet been investigated.

\section{The vehicle routing problem with synchronization constraints and support vehicle-dependent service times} 
\label{sec:problem}
In the following section, we present the assumptions made for the VRP with support vehicle-dependent service times (\ref{sec:Assumptions}), introduce two different mixed-integer linear programs (\ref{sec:MILP}), and also show how to transform the models into different variants (\ref{sec:Variants}). 

\subsection{Assumptions}
\label{sec:Assumptions}
We make the following assumption on the operations of primary and support vehicles:
\begin{itemize}
\item The fleets of primary and support vehicles are each homogeneous.
\item A primary vehicle requires at least one support vehicle to perform a service at a customer.
\item Primary and support vehicles travel on the same network with the same speed between locations.
\item Primary and support vehicles can visit a customer not more than once.
\item The operations of both vehicle types must be synchronized in terms of time and location. That is, service at a customer site cannot begin until the primary vehicle and all support vehicles assigned to that job have arrived. No support vehicle can join or leave an ongoing operation of a primary vehicle.
\item The number of support vehicles and service time are inversely proportional, i.e., more support vehicles lead to smaller service times due to higher utilization of the primary vehicle. There is a limit to the maximum number of support vehicles at a customer site due to the limited maximum utilization of a primary vehicle.
\end{itemize}

We minimize the latest completion time at all customer locations and refer to this time as makespan. This objective function considers two issues. First, time is crucial, for example, in harvesting to minimize the impact of unpredictable external circumstances, such as weather. Second, it implicitly maximizes the utilization of expensive primary vehicles, e.g., forage harvesters.

\subsection{Mathematical programming formulations}
\label{sec:MILP}
\subsubsection{Sets, parameters, decision variables}
\label{sec:SetParamsVariables}
The problem is defined on a complete directed Graph $G(N, A)$ with the set of nodes $N$ and the set of arcs $A$. The set of homogeneous primary vehicles is denoted by $K$, while $O$ represents the set of homogeneous support vehicles. Set $N$ is defined as $N = \left\{0\right\} \cup C \cup \left\{ n\right\}$ with node $0$ as start depot, node $n$ as end depot, and $C$ as the set of all customer nodes. To facilitate modeling and easily incorporate the prohibition of switching between primary vehicles for support vehicles, we generate the set $C$ by duplicating the set of original customer nodes $V$ for each primary vehicle $k \in K$. Therefore, each primary vehicle $k$ has its own set of customers $C_k$ and set of nodes $N_k = \left\{0\right\} \cup C_k \cup \left\{ n\right\}$. We define $N^+ = N \setminus \{ n \}$ ($N_k^+ = N_k \setminus \{n\}$) and $N^- = N \setminus \{0\}$ ($N_k^- = N_k \setminus \{0\}$) as sets where arcs may start or end.

We additionally introduce set $\bar{N}_j$ as the set of identical nodes of node $j \in C$ to keep track of which nodes correspond to the same original customer node. In addition, we use a similar set $\widetilde{N}_j$ as the set of nodes related to node $j$. We distinguish between $\bar{N}_j$ and $\widetilde{N}_j$ to integrate the ability to split the service at a node without making major changes to the model. In case split service is not allowed, exactly one node in $\bar{N}_j$ must be visited by a primary vehicle, and thus, we set $\widetilde{N}_j = \bar{N}_j$. If split service is allowed, all nodes in $\bar{N}_j$ may be visited and $\widetilde{N}_j = \left\{ j \right\}$.

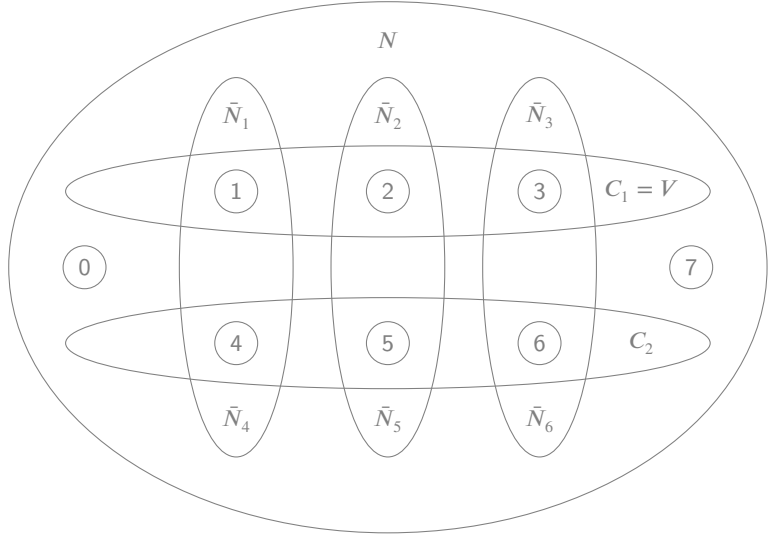
\begin{figure}
    \centering
    \begin{tikzpicture}[shorten >=1pt, auto, node distance=2cm]
    \tikzstyle{node_style} = [circle, draw]
    \tikzstyle{edge_style} = [draw=black, line width=1]
    \node[node_style] (v0) at (-4,0) {0};
    \node[node_style] (v1) at (-2,1) {1};
    \node[node_style] (v2) at (0,1) {2};
    \node[node_style] (v3) at (2,1) {3};
    \node[node_style] (v4) at (-2,-1) {4};
    \node[node_style] (v5) at (0,-1) {5};
    \node[node_style] (v6) at (2,-1) {6};
    \node[node_style] (v7) at (4,0) {7};
    \draw (0, 1) ellipse (4.25cm and 0.6cm);
    \draw (0,-1) ellipse (4.25cm and 0.6cm);
    \node (C) at (3.35, 1) {$C_1 = V$};
    \node (C) at (3.35,-1) {$C_2$};
    \draw (-2, 0) ellipse (0.75cm and 2.5cm);
    \node (N) at (-2, 2) {$\bar{N}_1$};
    \node (N) at (-2, -2) {$\bar{N}_4$};
    \draw (0, 0) ellipse (0.75cm and 2.5cm);
    \node (N) at (0, 2) {$\bar{N}_2$};
    \node (N) at (0, -2) {$\bar{N}_5$};
    \draw (2, 0) ellipse (0.75cm and 2.5cm);
    \node (N) at (2, 2) {$\bar{N}_3$};
    \node (N) at (2, -2) {$\bar{N}_6$};
    \draw (0, 0) ellipse (5cm and 3.5cm);
    \node (N) at (0, 3) {$N$};
\end{tikzpicture}
    \caption{Graph for an example with $|V| = 3$ original customer nodes and $|K| = 2$ primary vehicles.}
    \label{fig:example_graph}
\end{figure}

Consider Figure~\ref{fig:example_graph} for an example graph with original customer nodes $V = \left\{ 1,2,3 \right\}$, primary vehicles $K = \left\{ 1,2 \right\}$, and $n = 7$. The first three customer nodes $C_1 = \left\{ 1,2,3 \right\}$ belong to primary vehicle 1, whereas nodes $C_2 = \left\{ 4,5,6 \right\}$ belong to primary vehicle 2. Node 4 corresponds to the same original customer node as node 1. Therefore, $\bar{N}_1 = \left\{ 1,4 \right\} = \bar{N}_4$. If split service is not allowed, either node 1 or node 4 must be visited by a primary vehicle and $\widetilde{N}_1 = \widetilde{N}_4 = \bar{N}_1  = \bar{N_4}$. If split service is allowed, both nodes can be visited by a primary vehicle thus $\widetilde{N}_1 = \left\{ 1 \right\}$ and $\widetilde{N}_4 = \left\{ 4 \right\}$.

Each customer node $j \in C$ has a positive demand for service $d_j$ and a set of available operating modes $M_j = \left\{1,\hdots,b_j \right\}$ to fulfill this demand. The operating mode depends on the number of support vehicles that are used to serve the customer node. Demand $d_j$ is associated with the highest operating mode, i.e., the service requires $d_j$ time units if $b_j$ support vehicles are used. $p_j^m \in \left(0,1\right]$ is the productivity rate of mode $m \in M_j$ at customer node $j$ of the primary vehicle with $p_j^m = 1$ for the highest operating mode and corresponds to the primary vehicle utilization. Therefore, the real service time increases with fewer support vehicles due to decreasing productivity of the primary vehicle.

Each arc $\left( i,j \right)$ in arc set $A = \left\{ \left( i,j \right): i \in N^+, j \in N^-, j \notin \bar{N}_i \right\} $ has a non-negative travel time $\tau_{ij}$. We set all travel times to the end depot $n$ to zero, since we minimize the latest completion time at all customer nodes. Thus, the way back to depot is not included in the makespan determination. For the integer flow representation, we also introduce the arc capacity $\gamma_{ij}$ as the maximum number of support vehicles that are allowed to traverse arc $\left( i,j \right)$. Finally, $T$ represents the approximated maximum makespan used as Big M in some time-related constraints.

The decision variables to describe the two different mathematical programming formulations are given in Table~\ref{tab:decision_variables}. 

\begin{table}
	\begin{tabular}{lLl}
		\toprule
		\text{Model} & \text{Variable} & Description \\
		\midrule
		B, I &x_{ij} \in \left\{0,1\right\}& 1, if a primary vehicle travels directly from node $i$ to node $j$, \\
		     &                             & 0 otherwise \\
    		& q_j \in \left\{0,1 \right\} & 1, if a primary vehicle visits node $j$, 0 otherwise \\
    		& y_{j}^{m} \in  \left\{0,1\right\} & 1, if mode $m$ is used to serve node $j$, 0 otherwise\\
    		& t_j \in \mathbb{R} & latest time to start service at node $j$ \\
    		& s_j^m \in \mathbb{R} & service time at node $j$ in mode $m$ \\
    		\midrule
    		B & z_{ij}^o \in \left\{0,1\right\} & 1, if support vehicle $o$ travels directly from node $i$ to node $j$, \\
    		  &                                 & 0 otherwise \\
    		\midrule
    		I & v_{ij} \in \left\{0,1\right\} & 1, if any support vehicle travels directly from node $i$ to node $j$, \\ 
    		  &                               & 0 otherwise \\
    		& w_{ij} \in \mathbb{N} & number of support vehicles traveling directly from node $i$ to node $j$\\
		\bottomrule	
	\end{tabular}
	\caption{Decision variables used in the mathematical programming formulations with binary (B) and integer (I) flow  representation of support vehicles.} 
	\label{tab:decision_variables}
\end{table}

\subsubsection{Universal constraints}
Some constraints are the same in the binary and integer flow model. These mainly concern the routing of the primary vehicles and the service time at a node. We present them below and refer to them when introducing the models.
\begin{align}
	& \sum_{j \in C_k} x_{0j} \leq 1 \quad \forall~ k \in K \label{eq:outflow_primary_depot} \\
	& \sum_{i \in N_k^+} x_{ij} = q_j \quad \forall ~ k \in K, j \in C_k \label{eq:inflow_primary} \\
	& \sum_{i \in N_k^-} x_{ji} = q_j \quad \forall ~ k \in K, j \in C_k \label{eq:outflow_primary}\\
	& t_j \geq t_i + \sum_{h \in \widetilde{N}_i} \sum_{m \in M_h} s_h^m + \tau_{ij} \cdot x_{ij} - T \cdot \left( 1 - x_{ij} \right) \quad \forall ~ k \in K, i \in N_k^+, j \in N_k^-, i \neq j \label{eq:time_primary} \\
	& \sum_{i \in \bar{N}_j} \sum_{m \in M_i} s_i^m \cdot p_i^m = d_j \quad \forall ~ j \in V \label{eq:complete_demand_satisfaction} \\
	& \sum_{m \in M_j} y_j^m = q_j \quad \forall~ j \in C \label{eq:single_mode_node} \\
	& s_j^m \leq y_j^m \cdot \frac{d_j}{p_j^m} \quad \forall ~ j \in C, m \in M_j \label{eq:upper_bound_service} 
\end{align}
Constraints \eqref{eq:outflow_primary_depot} state that primary vehicle $k$ is allowed to enter the subset of its associated customer nodes $C_k$ at most once. Constraints \eqref{eq:inflow_primary} and \eqref{eq:outflow_primary} ensure that if primary vehicle $k$ visits customer node $j \in C_k$, it must leave this node. Constraints \eqref{eq:time_primary} set the latest start of service at node $j$, if primary vehicle $k$ travels directly from node $i \in N_k^+$ to node $j \in N_k^-$. Constraints \eqref{eq:complete_demand_satisfaction} guarantee that all demands are completely satisfied. Constraints \eqref{eq:single_mode_node} enforce that exactly one operating mode $m$ at customer node $j$ has to be used if and only if $j$ is visited by a primary vehicle. Constraints \eqref{eq:upper_bound_service} set an upper bound on the service time $s_j^m$ at customer node $j$ in mode $m$. The following constraints set the number of allowed visits for the model without and with split services.
\begin{subequations}
\label{eq:visits}
    \begin{align}
	& \sum_{i \in \bar{N}_j} q_i = 1 \quad \forall ~ j \in V \label{eq:all_customers_at_most_one_nosplit} \\
	& \sum_{i \in \bar{N}_j} q_i \geq 1 \quad \forall ~ j \in V \label{eq:all_customers_at_most_one_split}
	\end{align}
\end{subequations}
Without service splitting, exactly one node $i \in \bar{N}_j$ must be visited \eqref{eq:all_customers_at_most_one_nosplit}, while with service splitting, at least one node $i \in \bar{N}_j$ must be visited \eqref{eq:all_customers_at_most_one_split}.
\subsubsection{Binary flow representation of support vehicles}
Using the universal constraints \eqref{eq:outflow_primary_depot} - \eqref{eq:visits} introduced above, we can state the MILP with binary flow representation of support vehicles as follows:
\begin{align}
	& \min t_n \label{eq:objective_binary} \\
	& \text{s.t. } \eqref{eq:outflow_primary_depot} - \eqref{eq:visits} \notag \\
	& \sum_{m \in M_j} y_j^m \cdot m = \sum_{o \in O} \sum_{i \in N^+} z_{ij}^o \quad \forall ~ j \in C \label{eq:set_mode_binary} \\	
	& \sum_{j \in N^-} z_{0j}^o = 1 \quad \forall ~ o \in O  \label{eq:outflow_secondary_depot_binary}\\
	& \sum_{i \in N^+} z_{ij}^o = \sum_{i \in N^-} z_{ji}^o \quad \forall ~ j \in C, o \in O \label{eq:flow_conservation_binary} \\
	& t_j \geq t_i + \sum_{h \in \widetilde{N}_i} \sum_{m \in M_h} s_h^m + \tau_{ij} \cdot z_{ij}^o - T \cdot \left( 1 - z_{ij}^o \right) \quad \forall ~ o \in O, i \in N^+, j \in N^-, j \notin \bar{N}_i \label{eq:time_support_binary}
\end{align}
The objective function \eqref{eq:objective_binary} minimizes the latest completion time of all services at customer nodes. Besides universal constraints \eqref{eq:outflow_primary_depot} - \eqref{eq:upper_bound_service}, constraints \eqref{eq:set_mode_binary} define the operating mode that is used at node $j$ based on the number of support vehicles $o$ visiting node $j$. Constraints \eqref{eq:outflow_secondary_depot_binary} ensure that each support vehicle $o$ leaves the depot. Note that not all support vehicles need to be used, as a support vehicle $o$ may travel directly from node 0 to node $n$. Constraints \eqref{eq:flow_conservation_binary} are flow conservation constraints for each support vehicle $o$. Constraints \eqref{eq:time_support_binary} set the latest time to start the service at node $j$, if and only if support vehicle $o$ travels directly from node $i$ to node $j$. These constraints also synchronize the primary and support vehicles with respect to time, since we only use one variable $t_j$ for each node.

\subsubsection{Integer flow representation of support vehicles}
Similar to the binary flow representation, we use the universal constraints \eqref{eq:outflow_primary_depot} - \eqref{eq:visits} and can state the MILP with integer flow representation of support vehicles as follows: 
\begin{align}
	& \min t_n \label{eq:objective_integer}\\
	& \text{s.t. } \eqref{eq:outflow_primary_depot} - \eqref{eq:visits} \notag \\
	& \sum_{m \in M_j} y_j^m \cdot m = \sum_{i \in N^+} w_{ij} \quad \forall ~ j \in C \label{eq:set_mode_integer}\\	
	& \sum_{j \in N^-} w_{0j} = |O| \label{eq:outflow_secondary_depot_integer} \\
	& \sum_{i \in N^+} w_{ij} = \sum_{i \in N^-} w_{ji} \quad \forall ~ j \in C \label{eq:flow_conservation_integer} \\
	& w_{ij} \leq v_{ij} \cdot \gamma_{ij} \quad \forall ~ i \in N^+, j \in N^- \label{eq:arc_used_integer} \\
	& t_j \geq t_i + \sum_{h \in \widetilde{N}_i} \sum_{m \in M_h} s_h^m + \tau_{ij} \cdot v_{ij} - T \cdot \left( 1 - v_{ij} \right) \quad \forall ~ i \in N^+, j \in N^-, j \notin \bar{N}_i \label{eq:time_support_integer}
\end{align}
The objective function \eqref{eq:objective_integer} is the same as for the binary flow model and minimizes the latest completion time of all services at customer nodes. Constraints \eqref{eq:set_mode_integer} set the operating mode for service at node $j$ based on the number of support vehicles visiting node $j$. Constraints \eqref{eq:outflow_secondary_depot_integer} state that all support vehicles leave the depot. Similar to the binary flow model, not all support vehicles must be used at customer nodes. Constraints \eqref{eq:flow_conservation_integer} are flow conservation constraints for all customer nodes $j$. Therefore, the number of support vehicles entering node $j$ must equal the number of support vehicles leaving it. Constraints \eqref{eq:arc_used_integer} ensure that binary variables $v_{ij} = 1$ if a non-zero number of support vehicles travels on arc $(i,j)$. Furthermore, they restrict the number of support vehicles allowed to travel on arc $(i,j)$. Constraints \eqref{eq:time_support_integer} set the latest time to start the service at node $j$, if and only if any support vehicle travels directly from node $i$ to node $j$ and also ensure synchronization with respect to time between the two vehicle types.

\subsection{Model variants}
\label{sec:Variants}
The models introduced above can be easily transformed into several variants. Besides the binary (B) or integer (I) flow representation of support vehicles, these variants can represent different policies imposed by planners during the planning process and are based on:
\begin{itemize}
    \item switching of support vehicles between primary vehicles allowed (S) or not allowed (N)
    \item and splitting of services allowed (S) or not allowed (N).
\end{itemize}
A single variant is represented by a triple Flow\,|\,Switch\,|\,Split, e.g., integer representation, switching allowed, and splitting not allowed is denoted by \variant{I}{S}{N}. In our modeling approach, we can disable switching easily by enforcing 
\begin{equation}
    z^o_{ij} = 0 \quad \forall~ o \in O, k \in K, \{(i,j) | (i,j) \in A : i \in C_k \land j \notin C_k\}
\end{equation}
for the binary representation and by setting the arc capacity 
\begin{equation}
\gamma_{ij} = 0 \quad \forall~ k \in K, \{(i,j) | (i,j) \in A : i \in C_k \land j \notin C_k\}
\end{equation}
for the integer representation. Allowing split services can be achieved for both representations by defining the set of related nodes $\tilde{N}_j$ of customer $j \in C$ as $\tilde{N}_j = \{j\}$ as explained in Section~\ref{sec:SetParamsVariables} and by using constraints \eqref{eq:all_customers_at_most_one_split} instead of constraints \eqref{eq:all_customers_at_most_one_nosplit}. 

\section{Valid inequalities}
\label{sec:validInequalities}
The linear relaxations of the two mathematical programs are rather weak. We introduce several valid inequalities to strengthen the linear relaxation and improve run times and gaps at termination. We distinguish again between universal valid inequalities which can be used for both models and valid inequalities which are specific for the respective support vehicle representation. 
\subsection{Universal}
The following inequalities are valid for both support vehicle representations:
\begin{align}
    & t_n \geq \sum_{i \in N_k^+} \sum_{j \in N_k^-} \tau_{ij} \cdot x_{ij} + \sum_{j \in C_k} \sum_{m \in M_j} s_j^m \quad \forall ~ k \in K \
    	\label{eq:lower_bound_primary} \\
	&t_n \geq t_i + \sum_{h \in \widetilde{N}_i} \sum_{m \in M_h} s_h^m + \sum_{j \in N_k^-} \left(\tau_{ij} + \tau_{jn} \right) \cdot x_{ij} \quad \forall ~ k \in K, i \in N_k^+
	\label{eq:lower_bound_routing_primary} \\
	&t_j \geq \sum_{i \in N_k^+} \left(\tau_{0i} + \tau_{ij}\right) \cdot x_{ij} \quad \forall ~ k \in K, j \in N_k^-
	\label{eq:lower_bound_arrival_primary} \\
    & \sum_{j \in C_k}  x_{0j} \geq q_i \quad \forall k \in K, i \in C_k
	\label{eq:lower_bound_leaving_depot}
\end{align}
Inequalities \eqref{eq:lower_bound_primary} and \eqref{eq:lower_bound_routing_primary} set a lower bound on the makespan $t_n$. Inequalities \eqref{eq:lower_bound_primary} state that the makespan is greater than or equal to the sum of the travel and service times of each individual primary vehicle $k$. In contrast, inequalities \eqref{eq:lower_bound_routing_primary} are based on the latest departure time of primary vehicle $k$ from node $i$ and its earliest arrival at $n$, if a detour via any other node $j$ is considered. Similarly, inequalities \eqref{eq:lower_bound_arrival_primary} define lower bounds for the earliest arrival time at node $j$, if another node $i$ is visited by primary vehicle $k$ between the depot node $0$ and customer node $j$. Inequalities \eqref{eq:lower_bound_leaving_depot} state that primary vehicle $k$ must leave the depot, if $k$ visits any customer node $i$.

\subsection{Binary flow representation of support vehicles}
The following inequalities are valid for the binary flow representation:
\begin{align}
	t_n \geq \left( \sum_{o \in O} \sum_{i \in N^+} \sum_{j \in N^-} z_{ij}^o \tau_{ij} + \sum_{j \in C} \sum_{m \in M_j} s_j^m \cdot m \right) \cdot \frac{1}{|O|}
	\label{eq:binary_secondary_vehicle_lower_bound} \\
    \sum_{o \in O} \sum_{h \in N^+} \sum_{i \in \widetilde{N}_j} z_{hi}^o \leq \min\left\{|O|, b_j \right\} \quad \forall ~ j \in C
    \label{eq:binary_secondary_vehicle_upper_bound}
\end{align}

Inequalities \eqref{eq:binary_secondary_vehicle_lower_bound} set a lower bound on the makespan based on the time that all support vehicles spend traveling between nodes or serving customers. Unlike inequalities \eqref{eq:lower_bound_primary}, we cannot consider each support vehicle individually, since service times $s_j^m$ only apply to the primary vehicles. Inequalities (\ref{eq:binary_secondary_vehicle_upper_bound}) set an upper bound for the number of support vehicles leaving node $j$.

\subsection{Integer flow representation of support vehicles}
The following inequalities are valid for the integer flow representation:
\begin{align}
	t_n \geq \left( \sum_{i \in N^+} \sum_{j \in N^-} w_{ij} \tau_{ij} + \sum_{j \in C} \sum_{m \in M_j} s_j^m \cdot m \right) \cdot \frac{1}{|O|}
	\label{eq:integer_secondary_vehicle_lower_bound} \\
	\sum_{h \in N^+} \sum_{i \in \widetilde{N}_j} w_{hi} \leq \min\left\{|O|, b_j \right\} \quad \forall ~ j \in V
	\label{eq:integer_secondary_vehicle_upper_bound} \\
	w_{ij} \geq v_{ij} \quad \forall ~ i \in N^+, j \in N^-
	\label{eq:minimum_flow_on_arc}
\end{align}

Inequalities \eqref{eq:integer_secondary_vehicle_lower_bound} and \eqref{eq:integer_secondary_vehicle_upper_bound} have the same meaning as inequalities \eqref{eq:binary_secondary_vehicle_lower_bound} and \eqref{eq:binary_secondary_vehicle_upper_bound}. Inequalities (\ref{eq:minimum_flow_on_arc}) strengthen the relationship between binary and integer routing variables of support vehicles.

\section{Computational Experiments} \label{sec:ComputationalExperiments}

In this section, we first describe the creation of our test instances in Section~\ref{sec:testInstances}, then present the results of our experiments in Section~\ref{sec:results}. The latter includes the comparison of binary and integer flow representation (\ref{sec:ComparisonFlowRepresentation}), the influence of the number of primary and support vehicles on the makespan (\ref{sec:influenceMakespan}), and the analysis of different assumption and the resulting variants (\ref{sec:comparisonVariants}).

\subsection{Test instances and setting}
\label{sec:testInstances}
We generate instances with 5, 10, and 15 customer nodes for our computational experiments. All nodes are randomly distributed on a $100\times100$ plane. Travel times correspond to the Euclidean distances between nodes. The maximum number of support vehicles per primary vehicle at a customer node ranges from two to four. The demand of a node, i.e., the service time at a customer location with the maximum number of support vehicles, ranges from 20 to 50. 

For each number of customer nodes we vary the number of primary and support vehicles as follows:
\begin{equation*}
|K|=2 : |O|=\{4, 5, 6, 7\},\ |K|=3 : |O|=\{6, 7, 8, 9\}, \ |K|=4 : |O|=\{8, 9, 10, 11\}.
\end{equation*}
Thus, on average, there are at least two support vehicles per primary vehicle. In the remainder of this paper, we abbreviate a configuration consisting of number of customer nodes, number of primary vehicles, and number of support vehicles as |C|-|K|-|O|. With five instances per configuration, we consider 180 instances in total.

The MILP is implemented in Gurobi 9.1 via the Python 3.6 API. All instances are solved on an Intel(R) Xeon(R) CPU E5-2680 v3 with 2.50GHz, 8 cores, and 16 GB RAM. For each instance, the maximum runtime is two hours. Each instance is solved three times to account for the randomness of the solver.

The approximated makespan $T$ should be chosen as small as possible to strengthen the linear relaxations of the problem. The minimum value of $T$ is the optimal makespan. However, as this is initially unknown, we estimate $T$ by solving the simplest variant \variant{I}{N}{N} with a maximum runtime of 60 seconds. Note that this value of $T$ is valid for all other variants, and we do not include the runtime of this approximation in our results.

\subsection{Experiment results}
\label{sec:results}
We focus on aggregate data in the following sections to provide general insights. Detailed results of the experiments are given in Table~\ref{tab:experiment_results} in the Appendix.

\subsubsection{Comparison of binary and integer flow representation of support vehicles}
\label{sec:ComparisonFlowRepresentation}
In the following, we analyze and compare the binary and integer flow representation of support vehicles by applying the instances with five customer nodes. We use variants \variant{B}{S}{N} and \variant{I}{S}{N} (see Section~\ref{sec:Variants}) for these experiments. Table~\ref{tab:binaryinteger} shows a summary of the 180 runs (60 instances with three runs each) for both models. It contains the average of the objective function value (Makespan), the number of nonzero elements in the constraint matrix ({\#}Nonzeros), the number of branch-and-cut nodes explored ({\#}Nodes), the number of simplex iterations performed ({\#}Iterations), and the runtime to optimality (Time) in seconds.
\begin{table}
    \centering

    \begin{tabular}{lrrrrrrrr}
        \toprule
        Model  & $|V|$  & $|K|$ & Makespan & {\#}Nonzeros  &  {\#}Nodes & {\#}Iterations & Time\,[s]  \\

        \midrule
        \variant{B}{S}{N}         & 5 & 2 & 208.54 & 9870.00 & 187.97 &   5474.45 &  0.30 \\
                        &   & 3 & 156.98 & 30633.10 & 1248.10 &  54145.22 &  2.82  \\
                        &   & 4 & 132.29 & 74687.60 & 2332.85 &  136615.70 &  10.71 \\
        \midrule
        \variant{I}{S}{N}         & 5 & 2 & 208.54 & 4041.20 & 64.98 &     671.90 &  0.06  \\
                        &   & 3 & 156.98 & 8190.40 & 66.52 &     642.63 &  0.05 \\
                        &   & 4 & 132.29 & 14254.80 & 17.68 &     289.58 &  0.07 \\
    \bottomrule
    \end{tabular}
    \caption{Average of: makespan, number of nonzero elements in the constraint matrix, number of branch-and-cut nodes explored, total Simplex iterations, and computation time in seconds of models \variant{B}{S}{N} and \variant{I}{S}{N} for instances with five customer nodes.}
    \label{tab:binaryinteger}
\end{table}

The makespan is equal for all instances in both models. Thus, both models represent the same problem. However, significant differences can be observed in the other columns. The binary flow representation of support vehicles uses many more variables and constraints than the integer flow representation. This results in significantly larger problem sizes, represented by the number of nonzero elements in the constraint matrix. The problem size grows much faster in the binary representation as the number of primary vehicles increases. The smaller problem size in integer representation results in many fewer processed nodes, and fewer performed simplex iterations, ultimately leading to significantly shorter run times. The difference in runtime is particularly high for larger numbers of primary vehicles. Thus, the integer representation scales better with an increasing number of primary vehicles. Another reason for the longer runtimes is symmetries in the binary model due to the interchangeability of the support vehicles. However, a primary advantage of the binary model is the possibility of using support vehicles with different characteristics, resulting in a heterogeneous fleet.

The results of these experiments highlight the importance of an appropriate modeling approach. Here, many homogeneous support vehicles travel on a set of arcs, and it is common for multiple vehicles to travel on the same arcs. Describing this motion with integer variables has proven to be much more efficient. Therefore, we focus on the model with integer variables for the rest of the computational experiments.

\subsubsection{Influence of the number of primary and support vehicles on makespan}
\label{sec:influenceMakespan}

\begin{figure}
    \resizebox{\textwidth}{!}{\input{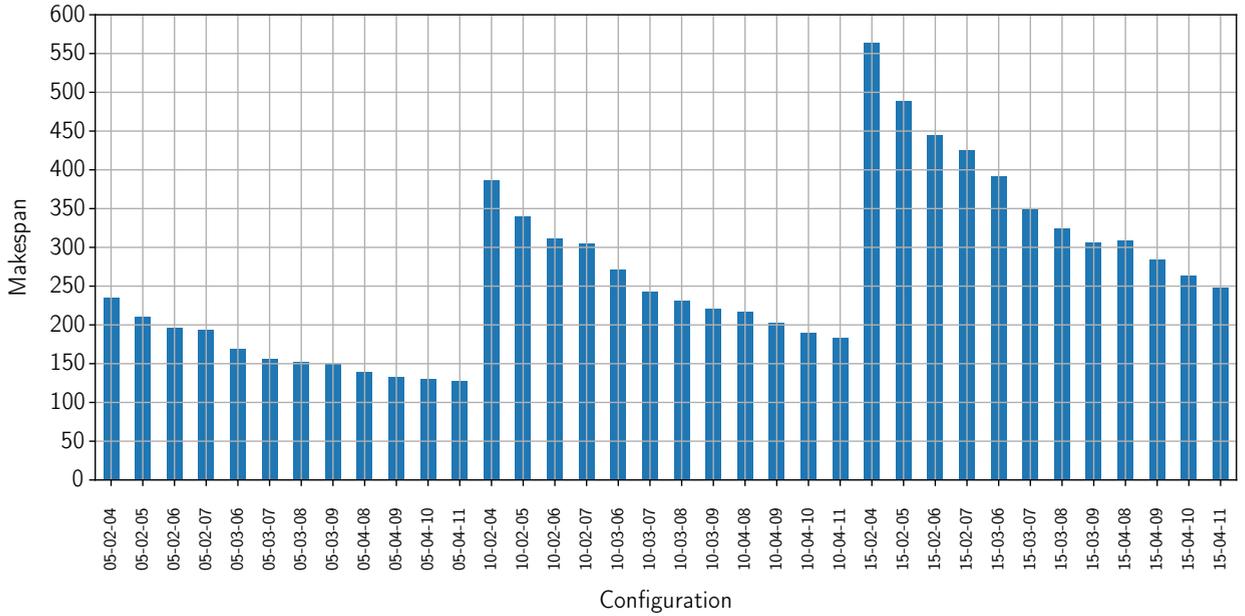}}
    \caption{Average makespan of variant \variant{I}{S}{N} for each configuration.}
    \label{fig:makespan}
\end{figure}

Figure~\ref{fig:makespan} shows the average makespan obtained with the model \variant{I}{S}{N} for each configuration. 
It can be deduced that increasing the number of vehicles, whether primary or support, reduces the makespan. However, the marginal benefit in terms of makespan reduction decreases for greater vehicle numbers. The extent of the reduction in the makespan also depends on the ratio of the two types of vehicles: Adding another vehicle of a relatively scarce resource results in a significantly greater makespan reduction.
For example, while the average makespan of configuration 05-02-04 can be reduced by 10.9\% by an additional support vehicle (05-02-05), the increase from six (05-02-06) to seven support vehicles (05-02-07) only results in a reduction of $1.1\%$. This behavior is observed due to a saturation effect when the average number of support vehicles per primary vehicle is greater than the maximum number required at some customer nodes for optimal primary vehicle utilization. Consequently, increasing the number of primary vehicles is more effective if the number of support vehicles is already high, e.g., adding another primary vehicle to configuration 05-03-08 reduces the makespan by $8.9\%$, while adding another primary vehicle to configuration 05-03-09 leads to a reduction of $11.5\%$.

Additionally, the ratio of primary vehicles to the number of customer nodes influences the makespan reduction. This behavior is likely due to vehicles interfering with each other on instances with a relatively small number of nodes. However, strong fluctuations in route lengths of the primary vehicles can be handled by corresponding adjustments in the number of support vehicles. In addition, if the number of nodes is larger, the average distance between nodes decreases.

\subsubsection{Analysis of model variants}
\label{sec:comparisonVariants}
\begin{figure}
    \captionsetup[subfigure]{aboveskip=1pt,belowskip=3pt}
    \centering
    \begin{subfigure}[c]{0.37\textwidth}
    \resizebox{!}{0.37\textheight}{\begin{tikzpicture}
\node[shape=circle,draw=black] (0) at (2.888492125933288,3.7641545598548394){0};
\node[shape=circle,draw=black] (1) at (4.166973884645603,2.7836574723837075){1};
\node[shape=circle,draw=black] (2) at (0.7544909863634021,4.97194166868202){2};
\node[shape=circle,draw=black] (3) at (3.6717431364592885,0.3169761664698412){3};
\node[shape=circle,draw=black] (4) at (2.454267225559507,1.286450484218962){4};
\node[shape=circle,draw=black] (5) at (0.10519823899346703,2.324794322362344){5};
\draw[dash pattern=on 5pt off 3pt, very thick] [->] (0) edge node {}(4);
\draw[dash pattern=on 4pt off 4pt, very thick] [->] (0) edge node {}(5);
\draw[dash pattern=on 2pt off 2pt, very thick] [->] (0) edge node {}(1);
\draw[dash pattern=on 5pt off 3pt, very thick] [->] (4) edge node {}(3);
\draw[dash pattern=on 2pt off 2pt, very thick] [->] (1) edge [bend right=30] node {}(2);
\draw [dash pattern=on 1pt off 1pt, thin] [->] (0) edge[bend right=30] node {} (4);
\draw [dash pattern=on 1pt off 1pt, thin] [->] (4) edge[bend right=30] node {} (3);
\draw [dash pattern=on 1pt off 2pt, thin] [->] (0) edge[bend right=20] node {} (4);
\draw [dash pattern=on 1pt off 2pt, thin] [->] (4) edge[bend right=20] node {} (3);
\draw [dash pattern=on 1pt off 3pt, thin] [->] (0) edge[bend right=20] node {} (5);
\draw [dash pattern=on 2pt off 1pt, thin] [->] (0) edge[bend right=40] node {} (1);
\draw [dash pattern=on 2pt off 1pt, thin] [->] (1) edge[bend right=50] node {} (2);
\draw [dash pattern=on 2pt off 2pt, thin] [->] (0) edge[bend right=30] node {} (1);
\draw [dash pattern=on 2pt off 2pt, thin] [->] (1) edge[bend right=40] node {} (2);
\draw [dash pattern=on 2pt off 3pt, thin] [->] (0) edge[bend right=20] node {} (1);
\node[text width=1cm] (1_0) at (-0.5,-0.4) {\scriptsize $k_1$};
\node[shape=rectangle,draw=black,minimum width=1.6973684210526316cm, minimum height=0.1cm] (1_4) at (2.106417166138044, -0.4){\tiny 4};
\node[text width=1cm] (2_0) at (-0.5,-0.8) {\scriptsize $k_2$};
\node[shape=rectangle,draw=black,minimum width=3.263157894736842cm, minimum height=0.1cm] (2_10) at (3.198301861156228, -0.8){\tiny 5};
\node[text width=1cm] (3_0) at (-0.5,-1.2000000000000002) {\scriptsize $k_3$};
\node[shape=rectangle,draw=black,minimum width=1.1578947368421053cm, minimum height=0.1cm] (3_11) at (1.384535685413595, -1.2000000000000002){\tiny 1};
\node[shape=rectangle,draw=black,minimum width=1.2105263157894737cm, minimum height=0.1cm] (1_3) at (4.3385239087229595, -0.4){\tiny 3};
\node[shape=rectangle,draw=black,minimum width=0.7894736842105263cm, minimum height=0.1cm] (3_12) at (4.385140420703859, -1.2000000000000002){\tiny 2};
\draw[dash pattern=on 5pt off 3pt, very thick] [->] (1_0) edge node {}(1_4);
\draw[dash pattern=on 4pt off 4pt, very thick] [->] (2_0) edge node {}(2_10);
\draw[dash pattern=on 2pt off 2pt, very thick] [->] (3_0) edge node {}(3_11);
\draw[dash pattern=on 5pt off 3pt, very thick] [->] (1_4) edge node {}(1_3);
\draw[dash pattern=on 2pt off 2pt, very thick] [->] (3_11) edge node {}(3_12);
\node[text width=1cm] (o1_0) at (-0.5,-1.6) {\scriptsize $o_1$};
\node[shape=rectangle,draw=black,minimum width=1.6973684210526316cm, minimum height=0.1cm] (o1_4) at (2.106417166138044, -1.6){\tiny 4};
\node[shape=rectangle,draw=black,minimum width=1.2105263157894737cm, minimum height=0.1cm] (o1_3) at (4.3385239087229595, -1.6){\tiny 3};
\node[text width=1cm] (o2_0) at (-0.5,-2.0) {\scriptsize $o_2$};
\node[shape=rectangle,draw=black,minimum width=1.6973684210526316cm, minimum height=0.1cm] (o2_4) at (2.106417166138044, -2.0){\tiny 4};
\node[shape=rectangle,draw=black,minimum width=1.2105263157894737cm, minimum height=0.1cm] (o2_3) at (4.3385239087229595, -2.0){\tiny 3};
\node[text width=1cm] (o3_0) at (-0.5,-2.4000000000000004) {\scriptsize $o_3$};
\node[shape=rectangle,draw=black,minimum width=3.263157894736842cm, minimum height=0.1cm] (o3_10) at (3.198301861156228, -2.4000000000000004){\tiny 5};
\node[text width=1cm] (o4_0) at (-0.5,-2.8000000000000003) {\scriptsize $o_4$};
\node[shape=rectangle,draw=black,minimum width=1.1578947368421053cm, minimum height=0.1cm] (o4_11) at (1.384535685413595, -2.8000000000000003){\tiny 1};
\node[shape=rectangle,draw=black,minimum width=0.7894736842105263cm, minimum height=0.1cm] (o4_12) at (4.385140420703859, -2.8000000000000003){\tiny 2};
\node[text width=1cm] (o5_0) at (-0.5,-3.2) {\scriptsize $o_5$};
\node[shape=rectangle,draw=black,minimum width=1.1578947368421053cm, minimum height=0.1cm] (o5_11) at (1.384535685413595, -3.2){\tiny 1};
\node[shape=rectangle,draw=black,minimum width=0.7894736842105263cm, minimum height=0.1cm] (o5_12) at (4.385140420703859, -3.2){\tiny 2};
\node[text width=1cm] (o6_0) at (-0.5,-3.6) {\scriptsize $o_6$};
\node[shape=rectangle,draw=black,minimum width=1.1578947368421053cm, minimum height=0.1cm] (o6_11) at (1.384535685413595, -3.6){\tiny 1};
\draw[dash pattern=on 1pt off 1pt, thin] [->] [->] (o1_0) edge node {}(o1_4);
\draw[dash pattern=on 1pt off 1pt, thin] [->] [->] (o1_4) edge node {}(o1_3);
\draw[dash pattern=on 1pt off 2pt, thin] [->] [->] (o2_0) edge node {}(o2_4);
\draw[dash pattern=on 1pt off 2pt, thin] [->] [->] (o2_4) edge node {}(o2_3);
\draw[dash pattern=on 1pt off 3pt, thin] [->] [->] (o3_0) edge node {}(o3_10);
\draw[dash pattern=on 2pt off 1pt, thin] [->] [->] (o4_0) edge node {}(o4_11);
\draw[dash pattern=on 2pt off 1pt, thin] [->] [->] (o4_11) edge node {}(o4_12);
\draw[dash pattern=on 2pt off 2pt, thin] [->] [->] (o5_0) edge node {}(o5_11);
\draw[dash pattern=on 2pt off 2pt, thin] [->] [->] (o5_11) edge node {}(o5_12);
\draw[dash pattern=on 2pt off 3pt, thin] [->] [->] (o6_0) edge node {}(o6_11);
\node[shape=rectangle,draw=white,minimum height=0.1cm] at (3.5, -3.8){};
\end{tikzpicture}}
    \subcaption{\variant{I}{N}{N} with makespan 187.9}
    \label{fig:INN}
    \end{subfigure}
    \begin{subfigure}[c]{0.37\textwidth}
    \resizebox{!}{0.37\textheight}{\begin{tikzpicture}
\node[shape=circle,draw=black] (0) at (2.888492125933288,3.7641545598548394){0};
\node[shape=circle,draw=black] (1) at (4.166973884645603,2.7836574723837075){1};
\node[shape=circle,draw=black] (2) at (0.7544909863634021,4.97194166868202){2};
\node[shape=circle,draw=black] (3) at (3.6717431364592885,0.3169761664698412){3};
\node[shape=circle,draw=black] (4) at (2.454267225559507,1.286450484218962){4};
\node[shape=circle,draw=black] (5) at (0.10519823899346703,2.324794322362344){5};
\draw[dash pattern=on 5pt off 3pt, very thick] [->] (0) edge[bend right=15] node {}(1);
\draw[dash pattern=on 4pt off 4pt, very thick] [->] (0) edge node {}(4);
\draw[dash pattern=on 2pt off 2pt, very thick] [->] (0) edge[bend left=10] node {}(1);
\draw[dash pattern=on 5pt off 3pt, very thick] [->] (1) edge node {}(3);
\draw[dash pattern=on 4pt off 4pt, very thick] [->] (4) edge node {}(5);
\draw[dash pattern=on 2pt off 2pt, very thick] [->] (1) edge[bend right=40] node {}(2);
\draw [dash pattern=on 1pt off 1pt, thin] [->] (0) edge[bend right=40] node {} (1);
\draw [dash pattern=on 1pt off 1pt, thin] [->] (1) edge[bend right=30] node {} (3);
\draw [dash pattern=on 1pt off 2pt, thin] [->] (0) edge[bend right=30] node {} (1);
\draw [dash pattern=on 1pt off 2pt, thin] [->] (1) edge[bend right=20] node {} (3);
\draw [dash pattern=on 1pt off 3pt, thin] [->] (0) edge[bend right=40] node {} (4);
\draw [dash pattern=on 1pt off 3pt, thin] [->] (4) edge[bend right=40] node {} (5);
\draw [dash pattern=on 2pt off 1pt, thin] [->] (0) edge[bend right=30] node {} (4);
\draw [dash pattern=on 2pt off 1pt, thin] [->] (4) edge[bend right=30] node {} (5);
\draw [dash pattern=on 2pt off 2pt, thin] [->] (0) edge[bend right=20] node {} (4);
\draw [dash pattern=on 2pt off 2pt, thin] [->] (4) edge[bend right=20] node {} (5);
\draw [dash pattern=on 2pt off 3pt, thin] [->] (0) edge[bend left=25] node {} (1);
\draw [dash pattern=on 2pt off 3pt, thin] [->] (1) edge[bend right=50] node {} (2);
\node[text width=1cm] (1_0) at (-0.5,-0.4) {\scriptsize $k_1$};
\node[shape=rectangle,draw=black,minimum width=1.537024687999259cm, minimum height=0.1cm] (1_1) at (1.5741006609921688, -0.4){\tiny 1};
\node[text width=1cm] (2_0) at (-0.5,-0.8) {\scriptsize $k_2$};
\node[shape=rectangle,draw=black,minimum width=1.131578947368421cm, minimum height=0.1cm] (2_9) at (1.823522429295937, -0.8){\tiny 4};
\node[text width=1cm] (3_0) at (-0.5,-1.2000000000000002) {\scriptsize $k_3$};
\node[shape=rectangle,draw=black,minimum width=0.39963483452779774cm, minimum height=0.1cm] (3_11) at (1.0054057342564398, -1.2000000000000002){\tiny 1};
\node[shape=rectangle,draw=black,minimum width=1.2105263157894737cm, minimum height=0.1cm] (1_3) at (4.205827886810601, -0.4){\tiny 3};
\node[shape=rectangle,draw=black,minimum width=1.087719298245614cm, minimum height=0.1cm] (2_10) at (4.217333051862291, -0.8){\tiny 5};
\node[shape=rectangle,draw=black,minimum width=1.5789473684210527cm, minimum height=0.1cm] (3_12) at (4.021617360494813, -1.2000000000000002){\tiny 2};
\draw[dash pattern=on 5pt off 3pt, very thick] [->] (1_0) edge node {}(1_1);
\draw[dash pattern=on 4pt off 4pt, very thick] [->] (2_0) edge node {}(2_9);
\draw[dash pattern=on 2pt off 2pt, very thick] [->] (3_0) edge node {}(3_11);
\draw[dash pattern=on 5pt off 3pt, very thick] [->] (1_1) edge node {}(1_3);
\draw[dash pattern=on 4pt off 4pt, very thick] [->] (2_9) edge node {}(2_10);
\draw[dash pattern=on 2pt off 2pt, very thick] [->] (3_11) edge node {}(3_12);
\node[text width=1cm] (o1_0) at (-0.5,-1.6) {\scriptsize $o_1$};
\node[shape=rectangle,draw=black,minimum width=1.537024687999259cm, minimum height=0.1cm] (o1_1) at (1.5741006609921688, -1.6){\tiny 1};
\node[shape=rectangle,draw=black,minimum width=1.2105263157894737cm, minimum height=0.1cm] (o1_3) at (4.205827886810601, -1.6){\tiny 3};
\node[text width=1cm] (o2_0) at (-0.5,-2.0) {\scriptsize $o_2$};
\node[shape=rectangle,draw=black,minimum width=1.537024687999259cm, minimum height=0.1cm] (o2_1) at (1.5741006609921688, -2.0){\tiny 1};
\node[shape=rectangle,draw=black,minimum width=1.2105263157894737cm, minimum height=0.1cm] (o2_3) at (4.205827886810601, -2.0){\tiny 3};
\node[text width=1cm] (o3_0) at (-0.5,-2.4000000000000004) {\scriptsize $o_3$};
\node[shape=rectangle,draw=black,minimum width=1.131578947368421cm, minimum height=0.1cm] (o3_9) at (1.823522429295937, -2.4000000000000004){\tiny 4};
\node[shape=rectangle,draw=black,minimum width=1.087719298245614cm, minimum height=0.1cm] (o3_10) at (4.217333051862291, -2.4000000000000004){\tiny 5};
\node[text width=1cm] (o4_0) at (-0.5,-2.8000000000000003) {\scriptsize $o_4$};
\node[shape=rectangle,draw=black,minimum width=1.131578947368421cm, minimum height=0.1cm] (o4_9) at (1.823522429295937, -2.8000000000000003){\tiny 4};
\node[shape=rectangle,draw=black,minimum width=1.087719298245614cm, minimum height=0.1cm] (o4_10) at (4.217333051862291, -2.8000000000000003){\tiny 5};
\node[text width=1cm] (o5_0) at (-0.5,-3.2) {\scriptsize $o_5$};
\node[shape=rectangle,draw=black,minimum width=1.131578947368421cm, minimum height=0.1cm] (o5_9) at (1.823522429295937, -3.2){\tiny 4};
\node[shape=rectangle,draw=black,minimum width=1.087719298245614cm, minimum height=0.1cm] (o5_10) at (4.217333051862291, -3.2){\tiny 5};
\node[text width=1cm] (o6_0) at (-0.5,-3.6) {\scriptsize $o_6$};
\node[shape=rectangle,draw=black,minimum width=0.39963483452779774cm, minimum height=0.1cm] (o6_11) at (1.0054057342564398, -3.6){\tiny 1};
\node[shape=rectangle,draw=black,minimum width=1.5789473684210527cm, minimum height=0.1cm] (o6_12) at (4.021617360494813, -3.6){\tiny 2};
\draw[dash pattern=on 1pt off 1pt, thin] [->] [->] (o1_0) edge node {}(o1_1);
\draw[dash pattern=on 1pt off 1pt, thin] [->] [->] (o1_1) edge node {}(o1_3);
\draw[dash pattern=on 1pt off 2pt, thin] [->] [->] (o2_0) edge node {}(o2_1);
\draw[dash pattern=on 1pt off 2pt, thin] [->] [->] (o2_1) edge node {}(o2_3);
\draw[dash pattern=on 1pt off 3pt, thin] [->] [->] (o3_0) edge node {}(o3_9);
\draw[dash pattern=on 1pt off 3pt, thin] [->] [->] (o3_9) edge node {}(o3_10);
\draw[dash pattern=on 2pt off 1pt, thin] [->] [->] (o4_0) edge node {}(o4_9);
\draw[dash pattern=on 2pt off 1pt, thin] [->] [->] (o4_9) edge node {}(o4_10);
\draw[dash pattern=on 2pt off 2pt, thin] [->] [->] (o5_0) edge node {}(o5_9);
\draw[dash pattern=on 2pt off 2pt, thin] [->] [->] (o5_9) edge node {}(o5_10);
\draw[dash pattern=on 2pt off 3pt, thin] [->] [->] (o6_0) edge node {}(o6_11);
\draw[dash pattern=on 2pt off 3pt, thin] [->] [->] (o6_11) edge node {}(o6_12);
\node[shape=rectangle,draw=white,minimum height=0.1cm] at (3.5, -3.8){};
\end{tikzpicture}}
    \subcaption{\variant{I}{N}{S} with makespan 182.8}
    \end{subfigure}
    \begin{subfigure}[c]{0.37\textwidth}
    \resizebox{!}{0.37\textheight}{\begin{tikzpicture}
\node[shape=circle,draw=black] (0) at (2.888492125933288,3.7641545598548394){0};
\node[shape=circle,draw=black] (1) at (4.166973884645603,2.7836574723837075){1};
\node[shape=circle,draw=black] (2) at (0.7544909863634021,4.97194166868202){2};
\node[shape=circle,draw=black] (3) at (3.6717431364592885,0.3169761664698412){3};
\node[shape=circle,draw=black] (4) at (2.454267225559507,1.286450484218962){4};
\node[shape=circle,draw=black] (5) at (0.10519823899346703,2.324794322362344){5};
\draw[dash pattern=on 5pt off 3pt, very thick] [->] (0) edge node {}(4);
\draw[dash pattern=on 4pt off 4pt, very thick] [->] (0) edge node {}(1);
\draw[dash pattern=on 2pt off 2pt, very thick] [->] (0) edge node {}(2);
\draw[dash pattern=on 5pt off 3pt, very thick] [->] (4) edge node {}(3);
\draw[dash pattern=on 2pt off 2pt, very thick] [->] (2) edge node {}(5);
\draw [dash pattern=on 1pt off 1pt, thin] [->] (0) edge[bend right=40] node {} (4);
\draw [dash pattern=on 1pt off 1pt, thin] [->] (4) edge[bend right=30] node {} (3);
\draw [dash pattern=on 1pt off 2pt, thin] [->] (0) edge[bend right=30] node {} (4);
\draw [dash pattern=on 1pt off 2pt, thin] [->] (4) edge[bend right=20] node {} (3);
\draw [dash pattern=on 1pt off 3pt, thin] [->] (0) edge[bend right=20] node {} (4);
\draw [dash pattern=on 1pt off 3pt, thin] [->] (4) edge[bend right=20] node {} (5);
\draw [dash pattern=on 2pt off 1pt, thin] [->] (0) edge[bend right=20] node {} (1);
\draw [dash pattern=on 2pt off 2pt, thin] [->] (0) edge[bend right=30] node {} (2);
\draw [dash pattern=on 2pt off 2pt, thin] [->] (2) edge[bend right=30] node {} (5);
\draw [dash pattern=on 2pt off 3pt, thin] [->] (0) edge[bend right=20] node {} (2);
\draw [dash pattern=on 2pt off 3pt, thin] [->] (2) edge[bend right=20] node {} (5);
\node[text width=1cm] (1_0) at (-0.5,-0.4) {\scriptsize $k_1$};
\node[shape=rectangle,draw=black,minimum width=1.131578947368421cm, minimum height=0.1cm] (1_4) at (1.8235224292959389, -0.4){\tiny 4};
\node[text width=1cm] (2_0) at (-0.5,-0.8) {\scriptsize $k_2$};
\node[shape=rectangle,draw=black,minimum width=3.473684210526316cm, minimum height=0.1cm] (2_6) at (2.5424304222556997, -0.8){\tiny 1};
\node[text width=1cm] (3_0) at (-0.5,-1.2000000000000002) {\scriptsize $k_3$};
\node[shape=rectangle,draw=black,minimum width=0.7894736842105263cm, minimum height=0.1cm] (3_12) at (1.6207782936817068, -1.2000000000000002){\tiny 2};
\node[shape=rectangle,draw=black,minimum width=1.2105263157894737cm, minimum height=0.1cm] (1_3) at (3.7727344350387493, -0.4){\tiny 3};
\node[shape=rectangle,draw=black,minimum width=1.087719298245614cm, minimum height=0.1cm] (3_15) at (4.2173330518623, -1.2000000000000002){\tiny 5};
\draw[dash pattern=on 5pt off 3pt, very thick] [->] (1_0) edge node {}(1_4);
\draw[dash pattern=on 4pt off 4pt, very thick] [->] (2_0) edge node {}(2_6);
\draw[dash pattern=on 2pt off 2pt, very thick] [->] (3_0) edge node {}(3_12);
\draw[dash pattern=on 5pt off 3pt, very thick] [->] (1_4) edge node {}(1_3);
\draw[dash pattern=on 2pt off 2pt, very thick] [->] (3_12) edge node {}(3_15);
\node[text width=1cm] (o1_0) at (-0.5,-1.6) {\scriptsize $o_1$};
\node[shape=rectangle,draw=black,minimum width=1.131578947368421cm, minimum height=0.1cm] (o1_4) at (1.8235224292959389, -1.6){\tiny 4};
\node[shape=rectangle,draw=black,minimum width=1.2105263157894737cm, minimum height=0.1cm] (o1_3) at (3.7727344350387493, -1.6){\tiny 3};
\node[text width=1cm] (o2_0) at (-0.5,-2.0) {\scriptsize $o_2$};
\node[shape=rectangle,draw=black,minimum width=1.131578947368421cm, minimum height=0.1cm] (o2_4) at (1.8235224292959389, -2.0){\tiny 4};
\node[shape=rectangle,draw=black,minimum width=1.2105263157894737cm, minimum height=0.1cm] (o2_3) at (3.7727344350387493, -2.0){\tiny 3};
\node[text width=1cm] (o3_0) at (-0.5,-2.4000000000000004) {\scriptsize $o_3$};
\node[shape=rectangle,draw=black,minimum width=1.131578947368421cm, minimum height=0.1cm] (o3_4) at (1.8235224292959389, -2.4000000000000004){\tiny 4};
\node[shape=rectangle,draw=black,minimum width=1.087719298245614cm, minimum height=0.1cm] (o3_15) at (4.2173330518623, -2.4000000000000004){\tiny 5};
\node[text width=1cm] (o4_0) at (-0.5,-2.8000000000000003) {\scriptsize $o_4$};
\node[shape=rectangle,draw=black,minimum width=3.473684210526316cm, minimum height=0.1cm] (o4_6) at (2.5424304222556997, -2.8000000000000003){\tiny 1};
\node[text width=1cm] (o5_0) at (-0.5,-3.2) {\scriptsize $o_5$};
\node[shape=rectangle,draw=black,minimum width=0.7894736842105263cm, minimum height=0.1cm] (o5_12) at (1.6207782936817068, -3.2){\tiny 2};
\node[shape=rectangle,draw=black,minimum width=1.087719298245614cm, minimum height=0.1cm] (o5_15) at (4.2173330518623, -3.2){\tiny 5};
\node[text width=1cm] (o6_0) at (-0.5,-3.6) {\scriptsize $o_6$};
\node[shape=rectangle,draw=black,minimum width=0.7894736842105263cm, minimum height=0.1cm] (o6_12) at (1.6207782936817068, -3.6){\tiny 2};
\node[shape=rectangle,draw=black,minimum width=1.087719298245614cm, minimum height=0.1cm] (o6_15) at (4.2173330518623, -3.6){\tiny 5};
\draw[dash pattern=on 1pt off 1pt, thin] [->] [->] (o1_0) edge node {}(o1_4);
\draw[dash pattern=on 1pt off 1pt, thin] [->] [->] (o1_4) edge node {}(o1_3);
\draw[dash pattern=on 1pt off 2pt, thin] [->] [->] (o2_0) edge node {}(o2_4);
\draw[dash pattern=on 1pt off 2pt, thin] [->] [->] (o2_4) edge node {}(o2_3);
\draw[dash pattern=on 1pt off 3pt, thin] [->] [->] (o3_0) edge node {}(o3_4);
\draw[dash pattern=on 1pt off 3pt, thin] [->] [->] (o3_4) edge node {}(o3_15);
\draw[dash pattern=on 2pt off 1pt, thin] [->] [->] (o4_0) edge node {}(o4_6);
\draw[dash pattern=on 2pt off 2pt, thin] [->] [->] (o5_0) edge node {}(o5_12);
\draw[dash pattern=on 2pt off 2pt, thin] [->] [->] (o5_12) edge node {}(o5_15);
\draw[dash pattern=on 2pt off 3pt, thin] [->] [->] (o6_0) edge node {}(o6_12);
\draw[dash pattern=on 2pt off 3pt, thin] [->] [->] (o6_12) edge node {}(o6_15);
\node[shape=rectangle,draw=white,minimum height=0.1cm] at (3.5, -3.8){};
\end{tikzpicture}}
    \subcaption{\variant{I}{S}{N} with makespan 180.9}
    \end{subfigure}
    \begin{subfigure}[c]{0.37\textwidth}
    \resizebox{!}{0.37\textheight}{\begin{tikzpicture}
\node[shape=circle,draw=black] (0) at (2.888492125933288,3.7641545598548394){0};
\node[shape=circle,draw=black] (1) at (4.166973884645603,2.7836574723837075){1};
\node[shape=circle,draw=black] (2) at (0.7544909863634021,4.97194166868202){2};
\node[shape=circle,draw=black] (3) at (3.6717431364592885,0.3169761664698412){3};
\node[shape=circle,draw=black] (4) at (2.454267225559507,1.286450484218962){4};
\node[shape=circle,draw=black] (5) at (0.10519823899346703,2.324794322362344){5};
\draw[dash pattern=on 5pt off 3pt, very thick] [->] (0) edge node {}(1);
\draw[dash pattern=on 4pt off 4pt, very thick] [->] (0) edge node {}(4);
\draw[dash pattern=on 2pt off 2pt, very thick] [->] (0) edge node {}(2);
\draw[dash pattern=on 5pt off 3pt, very thick] [->] (1) edge node {}(3);
\draw[dash pattern=on 4pt off 4pt, very thick] [->] (4) edge node {}(5);
\draw[dash pattern=on 2pt off 2pt, very thick] [->] (2) edge node {}(5);
\draw [dash pattern=on 1pt off 1pt, thin] [->] (0) edge[bend right=30] node {} (1);
\draw [dash pattern=on 1pt off 1pt, thin] [->] (1) edge[bend right=30] node {} (3);
\draw [dash pattern=on 1pt off 2pt, thin] [->] (0) edge[bend right=20] node {} (1);
\draw [dash pattern=on 1pt off 2pt, thin] [->] (1) edge[bend right=20] node {} (3);
\draw [dash pattern=on 1pt off 3pt, thin] [->] (0) edge[bend right=30] node {} (4);
\draw [dash pattern=on 1pt off 3pt, thin] [->] (4) edge[bend right=20] node {} (3);
\draw [dash pattern=on 2pt off 1pt, thin] [->] (0) edge[bend right=20] node {} (4);
\draw [dash pattern=on 2pt off 1pt, thin] [->] (4) edge[bend right=20] node {} (5);
\draw [dash pattern=on 2pt off 2pt, thin] [->] (0) edge[bend right=30] node {} (2);
\draw [dash pattern=on 2pt off 2pt, thin] [->] (2) edge[bend right=30] node {} (5);
\draw [dash pattern=on 2pt off 3pt, thin] [->] (0) edge[bend right=20] node {} (2);
\draw [dash pattern=on 2pt off 3pt, thin] [->] (2) edge[bend right=20] node {} (5);
\node[text width=1cm] (1_0) at (-0.5,-0.4) {\scriptsize $k_1$};
\node[shape=rectangle,draw=black,minimum width=1.736842105263158cm, minimum height=0.1cm] (1_1) at (1.674009369624121, -0.4){\tiny 1};
\node[text width=1cm] (2_0) at (-0.5,-0.8) {\scriptsize $k_2$};
\node[shape=rectangle,draw=black,minimum width=1.6973684210526316cm, minimum height=0.1cm] (2_9) at (2.1064171661380446, -0.8){\tiny 4};
\node[text width=1cm] (3_0) at (-0.5,-1.2000000000000002) {\scriptsize $k_3$};
\node[shape=rectangle,draw=black,minimum width=0.7894736842105263cm, minimum height=0.1cm] (3_12) at (1.6207782936817052, -1.2000000000000002){\tiny 2};
\node[shape=rectangle,draw=black,minimum width=0.8070175438596492cm, minimum height=0.1cm] (1_3) at (4.2038909181095905, -0.4){\tiny 3};
\node[shape=rectangle,draw=black,minimum width=0.5137587108035998cm, minimum height=0.1cm] (2_10) at (4.496142231825511, -0.8){\tiny 5};
\node[shape=rectangle,draw=black,minimum width=1.3746995919666212cm, minimum height=0.1cm] (3_15) at (4.065671791243999, -1.2000000000000002){\tiny 5};
\draw[dash pattern=on 5pt off 3pt, very thick] [->] (1_0) edge node {}(1_1);
\draw[dash pattern=on 4pt off 4pt, very thick] [->] (2_0) edge node {}(2_9);
\draw[dash pattern=on 2pt off 2pt, very thick] [->] (3_0) edge node {}(3_12);
\draw[dash pattern=on 5pt off 3pt, very thick] [->] (1_1) edge node {}(1_3);
\draw[dash pattern=on 4pt off 4pt, very thick] [->] (2_9) edge node {}(2_10);
\draw[dash pattern=on 2pt off 2pt, very thick] [->] (3_12) edge node {}(3_15);
\node[text width=1cm] (o1_0) at (-0.5,-1.6) {\scriptsize $o_1$};
\node[shape=rectangle,draw=black,minimum width=1.736842105263158cm, minimum height=0.1cm] (o1_1) at (1.674009369624121, -1.6){\tiny 1};
\node[shape=rectangle,draw=black,minimum width=0.8070175438596492cm, minimum height=0.1cm] (o1_3) at (4.2038909181095905, -1.6){\tiny 3};
\node[text width=1cm] (o2_0) at (-0.5,-2.0) {\scriptsize $o_2$};
\node[shape=rectangle,draw=black,minimum width=1.736842105263158cm, minimum height=0.1cm] (o2_1) at (1.674009369624121, -2.0){\tiny 1};
\node[shape=rectangle,draw=black,minimum width=0.8070175438596492cm, minimum height=0.1cm] (o2_3) at (4.2038909181095905, -2.0){\tiny 3};
\node[text width=1cm] (o3_0) at (-0.5,-2.4000000000000004) {\scriptsize $o_3$};
\node[shape=rectangle,draw=black,minimum width=1.6973684210526316cm, minimum height=0.1cm] (o3_9) at (2.1064171661380446, -2.4000000000000004){\tiny 4};
\node[shape=rectangle,draw=black,minimum width=0.8070175438596492cm, minimum height=0.1cm] (o3_3) at (4.2038909181095905, -2.4000000000000004){\tiny 3};
\node[text width=1cm] (o4_0) at (-0.5,-2.8000000000000003) {\scriptsize $o_4$};
\node[shape=rectangle,draw=black,minimum width=1.6973684210526316cm, minimum height=0.1cm] (o4_9) at (2.1064171661380446, -2.8000000000000003){\tiny 4};
\node[shape=rectangle,draw=black,minimum width=0.5137587108035998cm, minimum height=0.1cm] (o4_10) at (4.496142231825511, -2.8000000000000003){\tiny 5};
\node[text width=1cm] (o5_0) at (-0.5,-3.2) {\scriptsize $o_5$};
\node[shape=rectangle,draw=black,minimum width=0.7894736842105263cm, minimum height=0.1cm] (o5_12) at (1.6207782936817052, -3.2){\tiny 2};
\node[shape=rectangle,draw=black,minimum width=1.3746995919666212cm, minimum height=0.1cm] (o5_15) at (4.065671791243999, -3.2){\tiny 5};
\node[text width=1cm] (o6_0) at (-0.5,-3.6) {\scriptsize $o_6$};
\node[shape=rectangle,draw=black,minimum width=0.7894736842105263cm, minimum height=0.1cm] (o6_12) at (1.6207782936817052, -3.6){\tiny 2};
\node[shape=rectangle,draw=black,minimum width=1.3746995919666212cm, minimum height=0.1cm] (o6_15) at (4.065671791243999, -3.6){\tiny 5};
\draw[dash pattern=on 1pt off 1pt, thin] [->] [->] (o1_0) edge node {}(o1_1);
\draw[dash pattern=on 1pt off 1pt, thin] [->] [->] (o1_1) edge node {}(o1_3);
\draw[dash pattern=on 1pt off 2pt, thin] [->] [->] (o2_0) edge node {}(o2_1);
\draw[dash pattern=on 1pt off 2pt, thin] [->] [->] (o2_1) edge node {}(o2_3);
\draw[dash pattern=on 1pt off 3pt, thin] [->] [->] (o3_0) edge node {}(o3_9);
\draw[dash pattern=on 1pt off 3pt, thin] [->] [->] (o3_9) edge node {}(o3_3);
\draw[dash pattern=on 2pt off 1pt, thin] [->] [->] (o4_0) edge node {}(o4_9);
\draw[dash pattern=on 2pt off 1pt, thin] [->] [->] (o4_9) edge node {}(o4_10);
\draw[dash pattern=on 2pt off 2pt, thin] [->] [->] (o5_0) edge node {}(o5_12);
\draw[dash pattern=on 2pt off 2pt, thin] [->] [->] (o5_12) edge node {}(o5_15);
\draw[dash pattern=on 2pt off 3pt, thin] [->] [->] (o6_0) edge node {}(o6_12);
\draw[dash pattern=on 2pt off 3pt, thin] [->] [->] (o6_12) edge node {}(o6_15);
\node[shape=rectangle,draw=white,minimum height=0.1cm] at (3.5, -3.8){};
\end{tikzpicture}}
    \subcaption{\variant{I}{S}{S} with makespan 180.6}
    \end{subfigure}
    \caption{Example solution for an instance with five customer nodes, three primary and six support vehicles (05-03-06) for the four model variants.}
    \label{fig:graph_gantt_example_solution}
\end{figure}
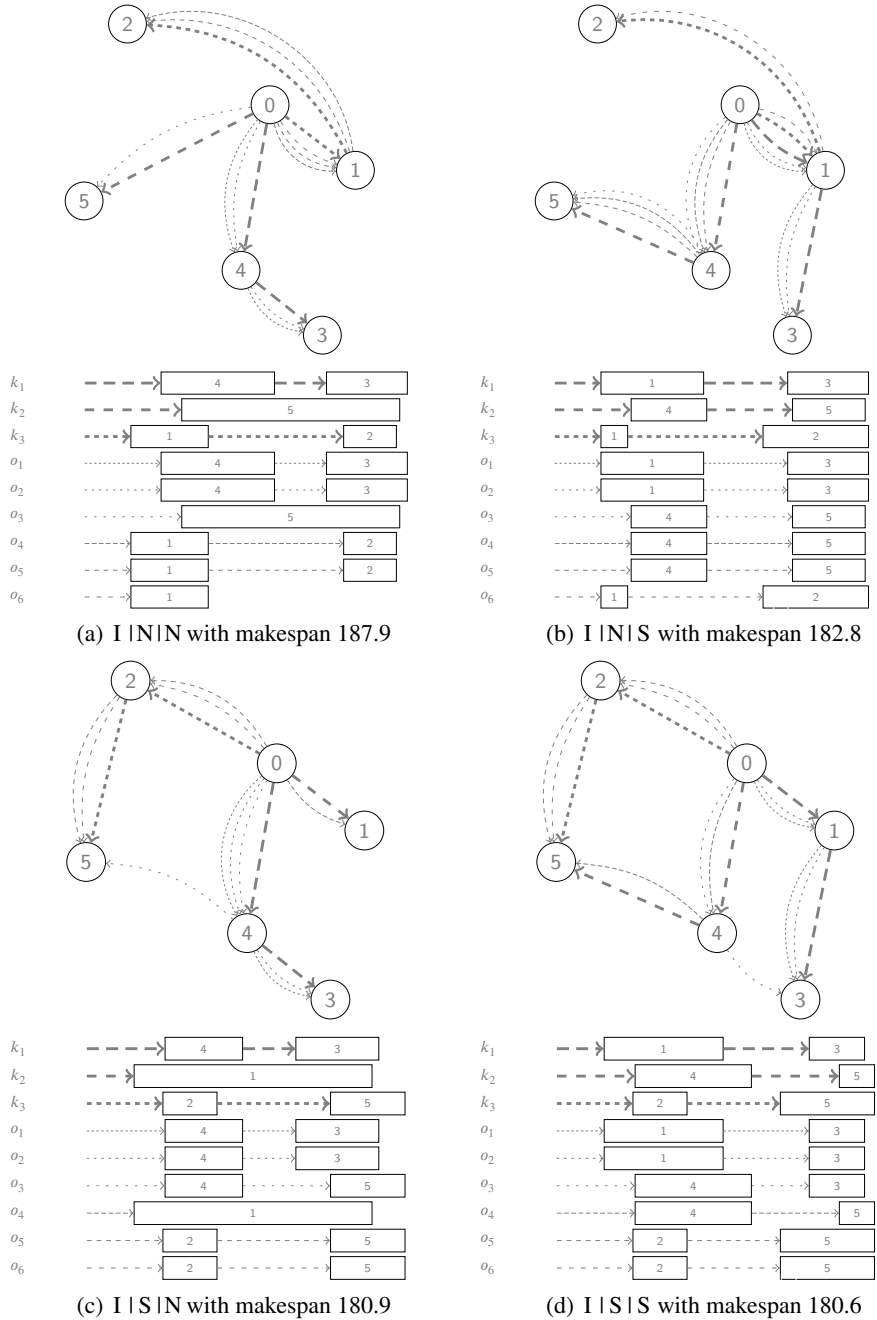

In addition to \variant{I}{S}{N}, we examine variants \variant{I}{N}{N}, \variant{I}{N}{S}, and \variant{I}{S}{S} to gain further insights into the problem characteristics and to provide additional managerial insights regarding the policies switching and splitting. 
Figure \ref{fig:graph_gantt_example_solution} depicts the optimal solutions of the four model variants for an instance with $V=\{1, 2, 3, 4, 5\}$, $K = \{k_1, k_2, k_3\}$ and $O = \{o_1, o_2, o_3, o_4, o_5, o_6\}$. Note that only the original customer nodes (set $V$) are displayed. The duplicated nodes to create the customer nodes of vehicles $k_2$ and $k_3$ (sets $C_{k_2}$ and $C_{k_3}$) are omitted in the figures since they overlap with the nodes in $V$. We show the solution on a graph and as a Gantt chart, highlighting the routing and scheduling aspects. 

All variants lead not only to different makespans but also to different solutions, e.g., regarding which customers are visited on a route or the number of support vehicles visiting a customer.
Variant \variant{I}{N}{N} results in an optimal makespan of $187.9$. Support vehicles cannot switch primary vehicles during makespan but can be distributed asymmetrically between primary vehicles. In this example, only support vehicle $o_3$ assists primary vehicle $k_3$ at customer 5, resulting in a relatively long service time at this node.
In \variant{I}{N}{S} the service at customer node 1 is split among primary vehicles $k_1$ and $k_3$. This results in a better distribution of services among the vehicles, and thus, the makespan can be reduced to $182.8$.
Allowing support vehicles to switch between primary vehicles (\variant{I}{S}{N}) results in a makespan of $180.9$, which in this example provides a larger reduction than allowing split services. While customer node 4 is served by primary vehicle $k_1$ with three support vehicles, $k_1$ continues its tour to customer node 3 with only two support vehicles. Support vehicle $o_3$ travels from customer node 4 to customer node 5, switching to the primary vehicle $k_3$.
The variant \variant{I}{S}{S} shares the same bottleneck as the previous model. However, with switching and splitting enabled, customer node 5 is served by primary vehicles $k_2$ and $k_3$. $k_3$ is assisted by two support vehicles, whereas $k_2$ is assisted by only one support vehicle, as support vehicle $o_3$ moves from node 4 to node 3 and switches to the primary vehicle $k_1$. As expected, this variant leads to the shortest makespan of $180.6$.

As shown in the example above, the policies regarding splitting and switching lead to different solutions and makespans. In the following, we analyze the general effects of switching and splitting. Table~\ref{tab:runtime_gap} presents the average computation time (Time) in seconds, the average optimality gap at termination (Gap) in percent, the number of optimally solved instances ($\#$Optimality), and the number of instances in which the variant achieves the best solution ($\#$Best Solution) for each model variant and problem size. Multiple variants can have the same best solution.
\begin{table}
\centering
\begin{tabular}{lrrrrr}
\toprule
Variant & $|V|$ &  Time\,[s] &    Gap\,[\%] & \#Optimality & \#Best Solution \\
\midrule
\variant{I}{N}{N}       & 5  &     0.04 &  0.00 & 60/60 & 15/60 \\
                        & 10 &     9.47 &  0.00 & 60/60 & 5/60 \\
                        & 15 &  1407.53 &  0.14 & 50/60 & 4/60 \\
\cmidrule[0.1pt](lr){3-6}
                        &   & 472.35   & 0.05    & 170/180 & 24/180 \\
\midrule
\variant{I}{N}{S}       & 5  &     0.32 &  0.00 & 60/60 & 50/60 \\
                        & 10 &    92.16 &  0.00 & 60/60 & 32/60 \\
                        & 15 &  3374.54 &  1.22 & 38/60 & 38/60 \\
\cmidrule[0.1pt](lr){3-6}
                        &   & 1155.67   & 0.41    & 158/180 & 120/180 \\
\midrule
\variant{I}{S}{N}       & 5  &     0.06 &  0.00 & 60/60 & 18/60 \\
                        & 10 &    35.14 &  0.00  & 60/60 & 12/60 \\
                        & 15 &  4397.69 &  3.49 & 22/60 & 13/6 \\
\cmidrule[0.1pt](lr){3-6}
                        &   & 1477.63   & 1.13    & 142/180 & 43/180 \\
\midrule
\variant{I}{S}{S}       & 5  &     2.55 &  0.00 & 60/60 & 60/60 \\
                        & 10 &  2054.55 &  1.14 & 43/60 & 57/60 \\
                        & 15 &  6092.41 &  4.77 & 12/60 & 39/60 \\
\cmidrule[0.1pt](lr){3-6}
                        &   & 2716.41   & 2.30    & 115/180 & 156/180 \\
\bottomrule
\end{tabular}
\caption{Average computation time, gap, number of instances optimally solved, number of best solutions in the comparison of variants for each model, and customer node number for a maximum runtime of 7200 seconds.}
\label{tab:runtime_gap}
\end{table}

The results in Table~\ref{tab:runtime_gap} highlight the trade-off between the different policies. On the one hand, prohibiting switching and splitting reduces the solution space and complexity of the problem and leads to problems that are easier to solve. Therefore, a high fraction of the instances can be solved to optimality by variant \variant{I}{N}{N}, while significantly fewer instances can be solved optimally by variant \variant{I}{S}{S}. Variants \variant{I}{N}{S} and \variant{I}{S}{N} are in between, where \variant{I}{N}{S} seems to be easier to solve than \variant{I}{S}{N} for larger instances. On the other hand, switching and splitting lead to better solutions with smaller makespans due to the increased solution space. Thus, \variant{I}{S}{S} most often provides the best solution. However, with 15 customer nodes, variants \variant{I}{N}{S} and \variant{I}{S}{S} achieve almost the same number of best solutions since \variant{I}{S}{S} often generates only suboptimal solutions in the given time limit. Therefore, reducing complexity in the planning process can be beneficial for larger instances, and allowing services to be split generally appears to be more effective than allowing switching between primary vehicles.

\begin{figure}
    \resizebox{\textwidth}{!}{\input{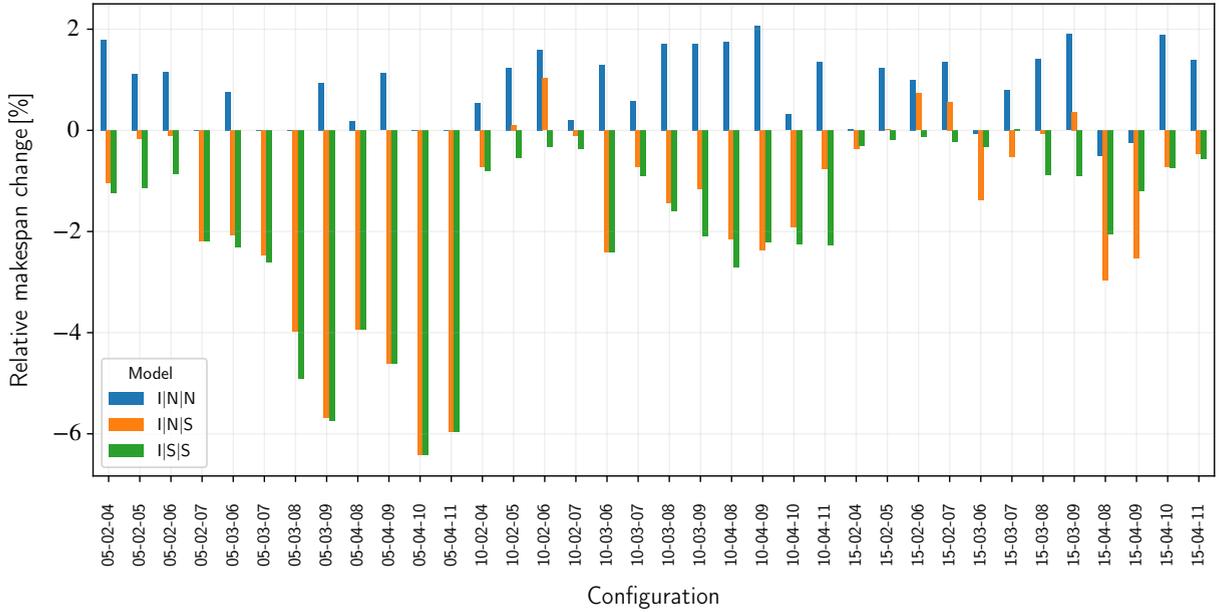}}
    \caption{Average relative makespan change in percent of variants \variant{I}{N}{N}, \variant{I}{N}{S}, and \variant{I}{S}{S} compared to \variant{I}{S}{N}.}
    \label{fig:relative_makespan}
\end{figure}

Figure~\ref{fig:relative_makespan} presents a closer look at the influence of the policies switching and splitting on the makespan. It shows the relative difference in the average makespan compared to variant \variant{I}{S}{N} for the remaining variants and each configuration of customer nodes, primary vehicles, and support vehicles.
Different effects can be observed for instances with five customer nodes, which can be solved to optimality by all variants. First, disabling switching between primary vehicles only slightly increases the makespan or does not increase the makespan at all. The latter usually occurs when more support vehicles are available. Thus, support vehicles are not a scarce resource that must be shared. Second, splitting customer services becomes more beneficial as the ratio of primary vehicles to customers increases. With splitting, the service at customer nodes can be more evenly distributed among the primary vehicles, while without splitting, the utilization of the primary vehicles is imbalanced. This behavior is particularly noticeable when many primary vehicles serve few customers.

This effect can also be seen to some extent for instances with ten and 15 customer nodes. However, results are less clear due to non-optimal solutions, especially with 15 customer nodes. The first effect of disabling switching cannot be observed with larger instances. Thus, the scarcity of support vehicles is determined by the number of primary vehicles and the number of customers per primary vehicle. With more customers per primary vehicle, switching allows for a better adjustment to the service times of customers than a fixed number of support vehicles per primary vehicle. This leads to a better distribution of services. In addition, customer locations are closer in larger instances, which can encourage switching between primary vehicles.

\section{Conclusions and future work} \label{sec:Conclusions}
This paper presents the VRP with synchronization constraints and support vehicle-dependent service times. The problem is inspired by the real-world planning process of forage harvesters and transporters, where both vehicle types must cooperate in harvesting a field. Therefore, this new problem considers the synchronized routing of two active resource types: primary (harvester) and support vehicles (transporter). Another critical aspect is that the number of support vehicles available for an operation affects the utilization rate of the primary vehicle. 

We introduce two novel MILPs, which differ in the representation of support vehicle flow. The underlying graph of the models enables a simple transformation into different model variants. The variants represent different policies that planners can impose during the planning process. They are based on switching support vehicles between primary vehicles and splitting customer services. We also introduce valid inequalities to strengthen the linear relaxation and improve the computation performance. 

Our computational experiments show that a representation of the support vehicle flow with integer decision variables is superior to a representation with binary decision variables. In addition, we examine the marginal benefit of support vehicles and primary vehicles. We show that the marginal benefit of support vehicles decreases with greater support vehicle numbers, whereby the marginal benefit of primary vehicles increases when support vehicles are added. Finally, we compare the different policy-based variants. We show that switching and splitting lead to problems that are more difficult to solve but result in better solutions with higher primary vehicle utilization. Furthermore, switching seems less beneficial than service splitting in our test instances.



For future work, a more general model may be developed in which support vehicles are enabled to leave or join an ongoing operation with reduced or increased primary vehicle utilization as a result, respectively. In addition, some practical applications will require solutions in a shorter computing time. This requires the development of corresponding heuristic or exact procedures to solve larger instances efficiently. Further research on the influence of the customer node density would supplement the managerial insights provided in this work.



\section*{Acknowledgements}

The work was partly carried out as a part of the IGF research project 19509 BR Landwirtschaftslogistik (Logistikdienstleistungen in der Landwirtschaft - Aufträge sicherer und ressourcenoptimal planen) and supported by Bundesvereinigung Logistik (BVL) e.V. The authors gratefully acknowledge the GWK support for funding this project by providing computing time through the Center for Information Services and HPC (ZIH) at TU Dresden.

\bibliographystyle{cas-model2-names}
\bibliography{References}   

\newpage
\begin{appendices}
\counterwithin{table}{section}

\section{Tables}

\centering
\begingroup
\small
\begin{longtable}{lrrrrrrr}
\toprule
Model & $|V|$ & $|K|$ & $|O|$ &  Makespan &  Time\,[s] &   Gap\,[\%] \\
\midrule
\endfirsthead 
\multicolumn{8}{c}{{\bfseries \tablename\ \thetable{} -- continued from previous page}}\\
\toprule
Model & $|V|$ & $|K|$ & $|O|$ &  Makespan &  Time\,[s] &   Gap\,[\%] \\
\midrule
\endhead 
\multicolumn{8}{r}{{Continued on next page}} \\
\midrule
\endfoot
\caption{Overview of the calculated average makespan, runtime, and gap for each configuration.}
\label{tab:experiment_results}
\endlastfoot
\variant{B}{S}{N} & 5  & 2 & 4  &    235.49 &     0.39 &   0.00  \\
                        &    &   & 5  &    209.87 &     0.34 &   0.00  \\
                        &    &   & 6  &    195.46 &     0.25 &   0.00  \\
                        &    &   & 7  &    193.33 &     0.22 &   0.00  \\
                        \cmidrule[0.1pt](lr){3-7}
                        &    & 3 & 6  &    169.55 &     4.21 &   0.00  \\
                        &    &   & 7  &    156.48 &     2.83 &   0.00  \\
                        &    &   & 8  &    152.22 &     2.56 &   0.00  \\
                        &    &   & 9  &    149.66 &     1.68 &   0.00  \\
                        \cmidrule[0.1pt](lr){3-7}
                        &    & 4 & 8  &    138.67 &    16.89 &   0.00  \\
                        &    &   & 9  &    132.49 &     9.35 &   0.00  \\
                        &    &   & 10 &    130.15 &     8.36 &   0.00  \\
                        &    &   & 11 &    127.84 &     8.23 &   0.00  \\
\midrule
\variant{I}{N}{N}       & 5  & 2 & 4  &    240.08 &     0.07 &   0.00  \\
                        &    &   & 5  &    211.76 &     0.07 &   0.00  \\
                        &    &   & 6  &    198.05 &     0.03 &   0.00  \\
                        &    &   & 7  &    193.33 &     0.02 &   0.00  \\
                        \cmidrule[0.1pt](lr){3-7}
                        &    & 3 & 6  &    170.94 &     0.06 &   0.00  \\
                        &    &   & 7  &    156.48 &     0.03 &   0.00  \\
                        &    &   & 8  &    152.22 &     0.02 &   0.00  \\
                        &    &   & 9  &    150.96 &     0.02 &   0.00  \\
                        \cmidrule[0.1pt](lr){3-7}
                        &    & 4 & 8  &    138.88 &     0.07 &   0.00  \\
                        &    &   & 9  &    133.92 &     0.05 &   0.00  \\
                        &    &   & 10 &    130.15 &     0.05 &   0.00  \\
                        &    &   & 11 &    127.84 &     0.05 &   0.00  \\
                        \cmidrule[0.1pt](lr){2-7}
                        & 10 & 2 & 4  &    388.33 &     2.82 &   0.00  \\
                        &    &   & 5  &    343.41 &     4.59 &   0.00  \\
                        &    &   & 6  &    316.69 &     0.65 &   0.00  \\
                        &    &   & 7  &    305.61 &     0.37 &   0.00  \\
                        \cmidrule[0.1pt](lr){3-7}
                        &    & 3 & 6  &    275.07 &    12.15 &   0.00  \\
                        &    &   & 7  &    244.52 &     4.30 &   0.00  \\
                        &    &   & 8  &    234.28 &     4.35 &   0.00  \\
                        &    &   & 9  &    225.03 &     1.11 &   0.00  \\
                        \cmidrule[0.1pt](lr){3-7}
                        &    & 4 & 8  &    220.32 &    59.16 &   0.00  \\
                        &    &   & 9  &    207.38 &    17.09 &   0.00  \\
                        &    &   & 10 &    190.66 &     4.76 &   0.00  \\
                        &    &   & 11 &    185.64 &     2.30 &   0.00  \\
                        \cmidrule[0.1pt](lr){2-7}
                        & 15 & 2 & 4  &    563.28 &   161.30 &   0.00  \\
                        &    &   & 5  &    493.63 &  1540.68 &   0.05  \\
                        &    &   & 6  &    449.30 &   288.91 &   0.00  \\
                        &    &   & 7  &    430.45 &     7.04 &   0.00  \\
                        \cmidrule[0.1pt](lr){3-7}
                        &    & 3 & 6  &    389.66 &  1352.40 &   0.00  \\
                        &    &   & 7  &    350.73 &   428.75 &   0.00  \\
                        &    &   & 8  &    328.27 &  1373.78 &   0.00  \\
                        &    &   & 9  &    312.42 &   236.89 &   0.00  \\
                        \cmidrule[0.1pt](lr){3-7}
                        &    & 4 & 8  &    306.10 &  4938.24 &   0.88  \\
                        &    &   & 9  &    281.25 &  2927.75 &   0.00  \\
                        &    &   & 10 &    266.87 &  2939.01 &   0.68  \\
                        &    &   & 11 &    251.87 &   695.56 &   0.00  \\
\midrule
\variant{I}{N}{S}              & 5  & 2 & 4  &    232.99 &     0.10 &   0.00  \\
                        &    &   & 5  &    209.11 &     0.09 &   0.00  \\
                        &    &   & 6  &    195.43 &     0.08 &   0.00  \\
                        &    &   & 7  &    189.30 &     0.06 &   0.00  \\
                        \cmidrule[0.1pt](lr){3-7}
                        &    & 3 & 6  &    166.11 &     0.29 &   0.00  \\
                        &    &   & 7  &    152.88 &     0.25 &   0.00  \\
                        &    &   & 8  &    146.43 &     0.24 &   0.00  \\
                        &    &   & 9  &    141.24 &     0.17 &   0.00  \\
                        \cmidrule[0.1pt](lr){3-7}
                        &    & 4 & 8  &    133.34 &     0.94 &   0.00  \\
                        &    &   & 9  &    126.64 &     0.71 &   0.00  \\
                        &    &   & 10 &    121.84 &     0.50 &   0.00  \\
                        &    &   & 11 &    120.26 &     0.43 &   0.00  \\
                        \cmidrule[0.1pt](lr){2-7}
                        & 10 & 2 & 4  &    383.47 &     5.45 &   0.00  \\
                        &    &   & 5  &    339.62 &     6.87 &   0.00  \\
                        &    &   & 6  &    314.90 &     1.47 &   0.00  \\
                        &    &   & 7  &    304.68 &     0.42 &   0.00  \\
                        \cmidrule[0.1pt](lr){3-7}
                        &    & 3 & 6  &    265.01 &    20.16 &   0.00  \\
                        &    &   & 7  &    241.29 &    32.54 &   0.00  \\
                        &    &   & 8  &    227.11 &    22.10 &   0.00  \\
                        &    &   & 9  &    218.57 &     8.73 &   0.00  \\
                        \cmidrule[0.1pt](lr){3-7}
                        &    & 4 & 8  &    211.89 &   470.16 &   0.00  \\
                        &    &   & 9  &    198.32 &   418.47 &   0.00  \\
                        &    &   & 10 &    186.37 &    65.93 &   0.00  \\
                        &    &   & 11 &    181.66 &    53.56 &   0.00  \\
                        \cmidrule[0.1pt](lr){2-7}
                        & 15 & 2 & 4  &    560.97 &   388.19 &   0.00  \\
                        &    &   & 5  &    487.84 &  1780.78 &   0.22  \\
                        &    &   & 6  &    448.09 &  2458.78 &   0.85  \\
                        &    &   & 7  &    427.04 &    15.72 &   0.00  \\
                        \cmidrule[0.1pt](lr){3-7}
                        &    & 3 & 6  &    384.54 &  2259.12 &   0.60  \\
                        &    &   & 7  &    346.14 &  2873.10 &   0.55  \\
                        &    &   & 8  &    323.42 &  3168.90 &   1.22  \\
                        &    &   & 9  &    307.72 &  4456.34 &   0.55  \\
                        \cmidrule[0.1pt](lr){3-7}
                        &    & 4 & 8  &    298.57 &  6094.96 &   2.54  \\
                        &    &   & 9  &    274.80 &  4795.19 &   2.39  \\
                        &    &   & 10 &    260.46 &  6102.47 &   3.41  \\
                        &    &   & 11 &    247.27 &  6100.99 &   2.31  \\
\midrule                        
\variant{I}{S}{N} & 5  & 2 & 4  &    235.49 &     0.09 &   0.00  \\
                        &    &   & 5  &    209.87 &     0.08 &   0.00  \\
                        &    &   & 6  &    195.46 &     0.03 &   0.00  \\
                        &    &   & 7  &    193.33 &     0.02 &   0.00  \\
                        \cmidrule[0.1pt](lr){3-7}
                        &    & 3 & 6  &    169.55 &     0.10 &   0.00  \\
                        &    &   & 7  &    156.48 &     0.05 &   0.00  \\
                        &    &   & 8  &    152.22 &     0.03 &   0.00  \\
                        &    &   & 9  &    149.66 &     0.03 &   0.00  \\
                        \cmidrule[0.1pt](lr){3-7}
                        &    & 4 & 8  &    138.67 &     0.10 &   0.00  \\
                        &    &   & 9  &    132.49 &     0.07 &   0.00  \\
                        &    &   & 10 &    130.15 &     0.05 &   0.00  \\
                        &    &   & 11 &    127.84 &     0.05 &   0.00  \\
                        \cmidrule[0.1pt](lr){2-7}
                        & 10 & 2 & 4  &    386.17 &    14.45 &   0.00  \\
                        &    &   & 5  &    339.25 &     6.30 &   0.00  \\
                        &    &   & 6  &    311.50 &     0.66 &   0.00  \\
                        &    &   & 7  &    304.98 &     0.38 &   0.00  \\
                        \cmidrule[0.1pt](lr){3-7}
                        &    & 3 & 6  &    271.45 &    96.90 &   0.00  \\
                        &    &   & 7  &    243.21 &    16.38 &   0.00  \\
                        &    &   & 8  &    230.44 &     8.96 &   0.00  \\
                        &    &   & 9  &    221.23 &     1.85 &   0.00  \\
                        \cmidrule[0.1pt](lr){3-7}
                        &    & 4 & 8  &    216.67 &   165.31 &   0.00  \\
                        &    &   & 9  &    203.18 &    90.25 &   0.00  \\
                        &    &   & 10 &    190.04 &    14.53 &   0.00  \\
                        &    &   & 11 &    183.16 &     5.73 &   0.00  \\
                        \cmidrule[0.1pt](lr){2-7}
                        & 15 & 2 & 4  &    563.35 &  4638.20 &   0.79  \\
                        &    &   & 5  &    488.01 &  4788.69 &   0.77  \\
                        &    &   & 6  &    444.81 &   589.68 &   0.00  \\
                        &    &   & 7  &    424.79 &     8.27 &   0.00  \\
                        \cmidrule[0.1pt](lr){3-7}
                        &    & 3 & 6  &    391.67 &  7200.05 &   6.14  \\
                        &    &   & 7  &    349.01 &  5853.17 &   4.46  \\
                        &    &   & 8  &    323.94 &  4026.11 &   0.90  \\
                        &    &   & 9  &    306.61 &   318.80 &   0.00  \\
                        \cmidrule[0.1pt](lr){3-7}
                        &    & 4 & 8  &    309.31 &  7200.03 &  11.54  \\
                        &    &   & 9  &    284.55 &  7200.06 &  10.85  \\
                        &    &   & 10 &    263.46 &  7200.06 &   6.42  \\
                        &    &   & 11 &    248.43 &  3749.17 &   0.03  \\
                        \midrule
\variant{I}{S}{S}              & 5  & 2 & 5  &    207.48 &     0.13 &   0.00  \\
                        &    &   & 6  &    193.84 &     0.08 &   0.00  \\
                        &    &   & 7  &    189.30 &     0.07 &   0.00  \\
                        \cmidrule[0.1pt](lr){3-7}
                        &    & 3 & 6  &    165.67 &     1.85 &   0.00  \\
                        &    &   & 7  &    152.62 &     1.20 &   0.00  \\
                        &    &   & 8  &    145.00 &     0.67 &   0.00  \\
                        &    &   & 9  &    141.13 &     0.27 &   0.00  \\
                        \cmidrule[0.1pt](lr){3-7}
                        &    & 4 & 8  &    133.34 &    17.80 &   0.00  \\
                        &    &   & 9  &    126.64 &     4.96 &   0.00  \\
                        &    &   & 10 &    121.84 &     1.94 &   0.00  \\
                        &    &   & 11 &    120.26 &     1.48 &   0.00  \\
                        \cmidrule[0.1pt](lr){2-7}
                        & 10 & 2 & 4  &    383.16 &   199.22 &   0.00  \\
                        &    &   & 5  &    337.40 &    58.22 &   0.00  \\
                        &    &   & 6  &    310.54 &     1.47 &   0.00  \\
                        &    &   & 7  &    303.88 &     0.42 &   0.00  \\
                        \cmidrule[0.1pt](lr){3-7}
                        &    & 3 & 6  &    265.18 &  3163.46 &   1.63  \\
                        &    &   & 7  &    241.05 &  1697.28 &   1.09  \\
                        &    &   & 8  &    226.75 &  1542.55 &   0.26  \\
                        &    &   & 9  &    216.48 &    82.61 &   0.00  \\
                        \cmidrule[0.1pt](lr){3-7}
                        &    & 4 & 8  &    211.42 &  6103.51 &   5.40  \\
                        &    &   & 9  &    198.83 &  6610.13 &   4.02  \\
                        &    &   & 10 &    185.70 &  2813.84 &   0.83  \\
                        &    &   & 11 &    178.95 &  2381.94 &   0.40  \\
                        \cmidrule[0.1pt](lr){2-7}
                        & 15 & 2 & 4  &    562.63 &  6618.11 &   1.60  \\
                        &    &   & 5  &    486.89 &  7140.60 &   2.28  \\
                        &    &   & 6  &    444.20 &  4324.58 &   1.18  \\
                        &    &   & 7  &    423.79 &   884.33 &   0.00  \\
                        \cmidrule[0.1pt](lr){3-7}
                        &    & 3 & 6  &    388.65 &  7200.03 &   6.25  \\
                        &    &   & 7  &    348.20 &  5936.23 &   5.58  \\
                        &    &   & 8  &    322.16 &  7001.60 &   4.94  \\
                        &    &   & 9  &    303.79 &  5203.24 &   1.54  \\
                        \cmidrule[0.1pt](lr){3-7}
                        &    & 4 & 8  &    302.37 &  7200.04 &   9.51  \\
                        &    &   & 9  &    279.56 &  7200.04 &   9.46  \\
                        &    &   & 10 &    260.83 &  7200.05 &   8.56  \\
                        &    &   & 11 &    247.56 &  7200.06 &   6.36  \\
\bottomrule
\end{longtable}
\endgroup

\end{appendices}

\end{document}